\documentclass[twoside,a4paper,11pt,limits]{amsart}
\usepackage{graphicx}
\usepackage{color}

\usepackage{enumitem}
\setenumerate{label={\rm (\alph{*})}}

\usepackage[rose,frak,hyperref,operator]{paper_diening}

\makeatletter
\long\def\unmarkedfootnote#1{{\long\def\@makefntext##1{##1}\footnotetext{#1}}}
\makeatother

\newtheorem{theorem}{Theorem}[section]
\newtheorem{lemma}[theorem]{Lemma}
\newtheorem{corollary}[theorem]{Corollary}

\newtheorem{remark}[theorem]{Remark}

\newtheorem{example}[theorem]{Example}

\numberwithin{equation}{section}
\providecommand{\meantmp}[2]{#1\langle{#2}#1\rangle}

\DeclareMathOperator{\rank}{rank}

\newcommand{\R}{\mathbb R}
\newcommand{\N}{\mathbb N}
\newcommand{\Div}{\divergence}

\newcommand{\rn}{\mathbb R^n}
\newcommand{\rmat}{\mathbb R^{n\times n}}

\newcommand{\dt}{\ensuremath{\,{\rm d} t}}

\newcommand{\dx}{\ensuremath{\,{\rm d} x}}

\DeclareMathOperator{\arsinh}{ar\,sinh}

\begin{document}

\title{Trace-free Korn inequalities in Orlicz spaces}

\begin{abstract}{Necessary and sufficient conditions are exhibited  for a  Korn-type inequality to hold between  (possibly different) Orlicz norms
 of the  gradient of vector-valued functions and
  of the deviatoric part of their symmetric gradients.
As a byproduct of our approach,   a positive answer is given to the
question  of the necessity of the same sufficient conditions in
related Korn-type inequalities for the full symmetric gradient, for
negative Orlicz-Sobolev norms, and for the gradient of the
\Bogovskii \, operator.
 }
\end{abstract}

\unmarkedfootnote {\par\noindent {\it Mathematics Subject
Classification: 46E35, 46E30.}
\par\noindent {\it Keywords:  Korn inequality, Orlicz spaces, symmetric gradient, trace-free symmetric gradient, singular
integrals.}
%
}


\author{D.~Breit, A.~Cianchi, L.~Diening}

\address{Department of Mathematics, Heriot-Watt University, Edinburgh EH14 4AS, UK}
\address{Dipartimento di Matematica e Informatica \lq\lq U. Dini", Universit\`a di Firenze, Viale Morgagni 67/A, 50134 Firenze, Italy}
\address{Mathematical Institute, University of Osnabr\"uck, Albrechtstr. 28a,
49076 Osnabr\"uck, Germany}

\email{d.breit@hw.ac.uk}
\email{cianchi@unifi.it}
\email{ldiening@uni-osnabrueck.de}

\maketitle
\section{Introduction}\label{intro}

The Korn inequality is a key tool in the analysis of mathematical
models for
  physical phenomena whose description only involves the symmetric part
   $\mathcal E {\bf u}$ of the distributional gradient $\nabla \bfu$ of
vector-valued functions ${\bf u}$. The theory of (generalized)
Newtonian fluids, and the classical theories of
plasticity and nonlinear elasticity
constitute paradigmatic examples in this connection.
 \par A plain form of the Korn inequality asserts that
 if $\Omega$ is an open bounded set in $\rn$, $n\geq 2$, and $1<p<\infty$,
then there exists a constant $C$ such that
\begin{equation}\label{intro1}
\int _\Omega |\nabla {\bf u}|^p\, dx  \leq C \int _\Omega | \mathcal
E {\bf u}|^p\, dx
\end{equation}
for every function $\bfu : \Omega \to \rn$ vanishing, in a suitable
sense, on $\partial \Omega$. Inequality \eqref{intro1} was
established by Korn in \cite{Korn} for $p=2$. Proofs of the general
case can be found in  \cite{Fu1, Go1,Go2, MM, Ne, Resh}. A fundamental reference for a
simple proof of the Korn's inequality in the modern setting is \cite{KO}.
Variants of inequality \eqref{intro1} are also available. For
instance, if $\Omega$ is connected and regular enough, a version of
\eqref{intro1} still holds if the boundary condition is dropped,
 and the left-hand side is
replaced with the ($p$-th power of the) distance, in the $L^p(\Omega
, \rmat)$ norm, of $\nabla {\bf u}$ from the space of skew-symmetric
matrices, namely the space of gradients of functions in the kernel
of the operator $\mathcal E$  \cite{DurMus04,DRS}. 
Let us incidentally mention that  nonlinear versions of Korn's inequality have been shown in \cite{FJM1,FJM2} (see also \cite{LMu} for further references).
An  extensive description of the historical background around Korn-type inequalities can be found in\cite{Ne}.
\par
The present paper is mainly concerned with  somewhat stronger,
closely related inequalities, where the symmetric gradient $\mathcal
E {\bf u}$ of a function $\bfu$ is replaced with its trace-free
part $\mathcal E^D {\bf u}$, also called deviatoric part of the
symmetric gradient. Inequalities of this kind  are critical  in the analysis of mathematical models for compressible fluids \cite{fei3,F} .
They also have important applications  to 
general relativity. Indeed, in the Cauchy formulation of the Einstein
gravitational field equations, the initial data have to satisfy the Einstein constraint
equations on a Riemannian manifold 
 of dimension $n > 2$  \cite{Ba}. One of these constraint equations amounts to the so-called momentum constraint equation, whose weak solutions can be obtained via
minimization of an energy functional depending on $\mathcal E^D {\bf u}$ \cite{Da}.
 Cosserat theory of elasticity is a further instance where trace-free Korn-type inequalities come into play \cite{FZ,Je,NeJe,NeJe2}.
\par A standard
trace-free Korn inequality reads
\begin{equation}\label{intro2}
\int _\Omega | \nabla {\bf u} |^p\, dx \leq C \int _\Omega |
\mathcal E^D {\bf u}|^p\, dx
\end{equation}
for every  $\bfu : \Omega \to \rn$ vanishing  on $\partial \Omega$.
A first proof of inequality \eqref{intro2}, and of its
  analogue  for functions with arbitrary boundary
values, can be found in  \cite{Resh}.
A comprehensive treatment of the  basic theory of the deviatoric Korn inequality (as well as of the standard Korn inequality) is offered, via modern analysis techniques, in
the monograph \cite{Resh2}. A simple
proof in case $p = 2$ was given by Dain \cite{Da}. We also refer to \cite{Sc2} for a proof in the case $1<p<\infty$. 
\\ A counterpart
of inequality \eqref{intro2} for functions with unprescribed
boundary values, where the left-hand side is replaced with the
$p$-th power of the distance from the space of gradients of
functions in the kernel of $\mathcal E^D$, takes a different form
depending on whether $n=2$ or $n \geq 3$. This kernel differs
substantially in the two cases, and, in particular, it agrees with
the whole space of holomorphic functions when $n=2$. The
inequalities in question require a distinct approach for $n=2$ and
for $n\geq 3$. In what follows, we shall focus on the case when $n
\geq 3$.
\par It is well known that
inequality \eqref{intro1}, and, a fortiori, inequality
\eqref{intro2} fail for the borderline values of the exponent $p$,
namely for $p=1$ \cite{Or} (see also \cite{BrD,CFM}) and $p=\infty$
(with integrals replaced with norms in $L^\infty (\Omega,\rmat)$)
\cite{BrD,LM}. The question thus arises of the validity of a version
of inequalities \eqref{intro1} and \eqref{intro2} where the role of
the power $t^p$ is played by a more general nonnegative convex
function $A(t)$ vanishing for $t=0$, briefly a Young function.  This
amounts to enlarging the class of $L^p$ norms of $\mathcal E \bfu$
and $\nabla \bfu$ in the Korn inequality to include the norms in the
Orlicz spaces $L^A(\Omega,\rmat)$.
  Korn-type inequalities in Orlicz spaces are relevant in the analysis of
mathematical models governed by strong  nonlinearities of
non-polynomial type.
\par  A  Korn inequality, for the symmetric gradient, in an
Orlicz space associated with a Young function $A$, is known to hold
if \cite{DRS,Fu2}, and only if \cite{BrD},  the Young function $A$
satisfies the so called $\Delta _2$ and $\nabla _2$ conditions near
infinity. Inequalities involving special Young functions fulfilling
the $\Delta _2$ and $\nabla _2$ conditions were earlier established
in \cite{AcMi} and \cite{BMM}. Loosely speaking, these conditions
amount to requiring that $A$ has a uniform rate of growth near
infinity, which is neither too slow, nor too rapid. On the other
hand, imposing the $\Delta _2$ and $\nabla _2$ conditions rules out
certain models in continuum mechanics. For instance,   the
nonlinearities appearing in the Prandt-Eyring fluids
\cite{BrDF,FuS}, and in models for plastic materials with
logarithmic hardening \cite{FrS}  are described  by a Young function
$A(t)$ that grows like $t\log(1+t)$ near infinity, and hence
violates the $\nabla_2$ condition. Young functions with fast growth,
which do not fulfil the $\Delta _2$ condition, are well suited to
model the behavior of fluids in certain liquid body armors
\cite{HRJ, SMB, W}.
\par
A  general  Orlicz version of the Korn inequality has been
established in \cite{Ci2}. In that paper,  it is shown that a Korn-type inequality for the symmetric gradient $\mathcal E \bfu$ in
$L^A(\Omega,\rmat)$ still holds, even if the $\Delta _2$ and $\nabla _2$
conditions on $A$ are dropped, provided that the norm of $\nabla
\bfu$ is taken in a possibly different Orlicz space $L^B(\Omega,\rmat)$.
The Young functions $A$ and $B$ have to be suitably balanced, in
such a way that the norm in  $L^B(\Omega,\rmat)$ turns out to be slightly
weaker than that in $L^A(\Omega,\rmat)$ when  $A$ does not fulfil either
the $\Delta _2$ condition, or the $\nabla _2$ condition near
infinity.
\par Here, we deal instead with trace-free Korn-type inequalities. Again, a priori arbitrary Orlicz spaces are allowed. In their
basic formulation for trial functions $\bfu$ vanishing on $\partial
\Omega$, the inequalities in question read
\begin{equation}\label{intro5}
\int _\Omega B\big(| \nabla {\bf u} |\big)\, dx \leq  \int _\Omega
A\big(C | \mathcal E^D {\bf u}|\big)\, dx.
\end{equation}
Their counterparts for functions $\bfu$ with unrestricted boundary
values, on a sufficiently regular connected and bounded open set
$\Omega$, take the form
\begin{equation}\label{intro6}
\inf_{\bfw\in\Sigma}\int_\Omega B(|\nabla\bfu-\nabla\bfw|)\,dx\leq
\int_\Omega A(C|\mathcal E^D\bfu|)\,dx\,,
\end{equation}
where $\Sigma$ denotes the kernel of the operator $\mathcal E^D$.
Our main result amounts to necessary and sufficient balance
conditions on the Young functions $A$ and $B$ for inequalities
\eqref{intro5} and \eqref{intro6} -- or slight variants of theirs
involving norms -- to hold.  It provides a comprehensive framework
for genuinely new trace-free Korn-type inequalities in borderline
customary and unconventional Orlicz spaces.
  Examples are exhibited in Section \ref{main} below. In particular, our characterization
 recovers
the fact that \eqref{intro5} and \eqref{intro6} hold with $B=A$ if
\cite{FB, BrThesis, BrSc}, and only if \cite{BrD} the function
$A$ fulfils both the $\Delta _2$ and the $\nabla _2$ condition near
infinity.
%

Let us emphasize that
the necessary and sufficient conditions for $A$ and $B$ to support
the Orlicz-Korn inequalities \eqref{intro5} and \eqref{intro6} turn
out to agree with those required in \cite{Ci2} for the Orlicz-Korn
inequalities for the standard symmetric gradient.
 In fact,  the necessity of the conditions for the former inequalities   follows via a proof of the
necessity of the same conditions  for the latter inequalities, an
issue which was left open in \cite{Ci2}. Another interesting
consequence is that we are now also in a position to derive the
necessity of parallel conditions on the Young functions appearing in
inequalities for negative Orlicz-Sobolev norms, and in inequalities
for the \Bogovskii \, operator in Orlicz spaces. These inequalities
have recently been established in \cite{BrCi} in connection with the
study of elliptic systems, with non-polynomial nonlinearities, in
 fluid mechanics.
\par To give an idea of the possible use of the results of this paper, we conclude this section with an outline of a model in fluid mechanics, for non-Newtonian fluids, where Korn inequalities, and trace-free Korn  inequalities in Orlicz spaces come into play.  
The stationary flow of an isentropic compressible fluid in a bounded domain $\Omega\subset\R^3$ can be described by the system
\begin{align}\label{eq:4.1}
\begin{cases}
-\Div \bfS+\Div\big(\varrho\bfu\otimes\bfu\big)+\nabla \pi=\varrho\,\bff \quad & \hbox{in $\Omega$,}\\
\Div (\varrho\bfu)=0 & \hbox{in $\Omega$,}\end{cases}
\end{align}
which accounts for  the balance of mass and momentum. Here, the velocity field $\bfu: \Omega\rightarrow\R^3$ and the density $\varrho:\Omega\rightarrow\R$ of the fluid are the unknown, whereas $\bff: \Omega\rightarrow\R^3$ is a given system of volume forces. The deviatoric stress tensor $\bfS: \Omega\rightarrow\R^{n\times n}$ and the pressure $\pi: \Omega\rightarrow\R$ have to be related to $\bfu$ and $\varrho$ by constitutive laws. A general model for non-Newtonian fluids takes the form
\begin{align*}
\bfS=\mu(|\mathcal E^D\bfu|)\mathcal E^D\bfu+\nu(|\Div\,\bfu|)(\Div\,\bfu)I,
\end{align*}
where $\mu,\nu:[0,\infty)\rightarrow[0,\infty)$ are given  functions, and  $I$ is the identity matrix -- see for instance \cite{MaRa} and \cite{FeMa}. If $\nu$ grows more slowly than $\mu$ (in fact, typically $\nu$ can even vanish), and the function $s\mu (s)$ is  non-decreasing, then  a priori estimates only imply that $\mathcal E^D\bfu\in L^A(\Omega)$, where
$$A(t)=\int_0^t s\,\mu(s) \,ds \quad \hbox{for $t \geq 0$}.$$
The  natural question that arises is to what extent the degree of integrability of $\mathcal E^D \bfu$ is inherited by  $\nabla u$.  This amounts to exibiting an optimal mutual dependence between the Young functions $A$ and $B$ in inequality  \eqref{intro5} or \eqref{intro6}.  Let us emphasize   that the mathematical literature about general non-Newtonian compressible fluids is quite limited. This is mainly due to the fact that an analogue to the existence theory from \cite{Li2} seems to be presently out of reach.\par
A much richer theory   is available in the incompressible case, corresponding to a constant density $\varrho$ in \eqref{eq:4.1}, and hence to the divergence-free constraint $\Div \bfu =0$. This implies that, $\mathcal E^D\bfu = \mathcal E\bfu$.
In the classical Prandtl-Eyring model, introduced in \cite{E}, the constitutive law reads 
\begin{equation}
  \label{0.3}\bfS=\eta_0 \frac{\arsinh (\lambda |\mathcal E\bfu|)}{\lambda |\mathcal E\bfu|}\mathcal E\bfu\,,
\end{equation}
where  $\eta_0$ and $\lambda$ are positive physical parameters. Since the function $\arsinh(t)$ behaves like $\log(1+t)$ near zero and infinity, the natural function space for the solutions $\bfu$ is obtained by requiring that $\mathcal E$ belongs to the Orlicz space $L\log L(\Omega,\R^{3\times 3})$. Theorem \ref{thm:plainkorn}, Section  \ref{main}, tells us that
\begin{align*}
\mathcal E\bfu\in L\log L(\Omega,\R^{3\times 3})\quad \hbox{implies that}\quad \nabla\bfu\in L^1(\Omega,\R^{3\times 3}),
\end{align*}
the space $L^1(\Omega,\R^{3\times 3})$ being optimal. This underlines the difficulties in the existence theory developed in \cite{BrDF} for stationary Prandtl-Eyring fluids, namely those satisfying   \eqref{eq:4.1} with $\varrho$ constant and $\bfS$ given by  \eqref{0.3}.
\\ The Bingham model amounts to the constitutive law
\begin{align}\label{eq:bing}
\bfS=\mu_0\mathcal E\bfu+\mu_\infty\frac{\mathcal E\bfu}{|\mathcal E\bfu|},
\end{align}
for some positive physical constants $\mu_0$ and $\mu_\infty$ -- see e.g. \cite{AsMa}. It is shown in \cite{FuS} that weak solutions
to the stationary Bingham model, consisting in \eqref{eq:4.1} with $\varrho$ constant and $\bfS$ obeying \eqref{eq:bing}, are such that $\mathcal E\bfu\in L^\infty (\Omega,\R^{3\times 3})$, at least locally. The resulting degree of inegrability is provided by Theorem \ref{thm:plainkorn}. It asserts that 
\begin{align*}
\mathcal E\bfu\in L^{\infty}(\Omega,\R^{3\times 3})\quad \hbox{implies that} \quad \nabla\bfu\in \exp L(\Omega,\R^{3\times 3})\,,
\end{align*}
and  the space $\exp L(\Omega,\R^{3\times 3})$ is optimal.

\section{Function spaces}\label{orlicz}

This section collects some definitions and basic results from the
theory of Orlicz and Orlicz-Sobolev spaces, as well as of their
versions for the symmetric, and trace-free symmetric gradient. For a
comprehensive treatment of the theory of Orlicz spaces we refer to
\cite{RR1,RR2}.
\par
 A function $A: [0, \infty) \to [0, \infty]$ is called a
Young function if it is convex, left-continuous, vanishing at $0$, and neither
identically equal to $0$, nor to $\infty$. Thus, with any such
function,  it is uniquely associated a (nontrivial) non-decreasing
left-continuous function $a:[0, \infty) \rightarrow [0, \infty]$
such that
\begin{equation}\label{B.4}
A(t) = \int _0^t a(r)\,dr \qquad {\rm for} \,\, t\geq 0.
\end{equation}
The Young conjugate $\widetilde{A}$ of $A$ is the Young function
defined by
$$\widetilde{A}(t) = \sup \{rt - A(r):\,r\geq 0\} \qquad {\rm for}\qquad  t\geq 0\,.$$
Note the representation formula
$$\widetilde A (t) = \int _0^t a^{-1}(r)\,dr \qquad {\rm for} \,\, t\geq
0,$$ where $a^{-1}$ denotes the (generalized) left-continuous
inverse of $a$. One has that
\begin{equation}\label{AAtilde}
r \leq A^{-1}(r)\widetilde A^{-1}(r) \leq 2 r \quad \hbox{for $r
\geq 0$,}
\end{equation}
where $A^{-1}$ denotes the (generalized) right-continuous inverse of
$A$.
Moreover,
\begin{equation}\label{B.7'}
\widetilde{\!\widetilde A\,}=A\,
\end{equation}
for any Young function $A$.
 If $A$ is any Young function and $\lambda \geq 1$, then
\begin{equation}\label{lambdaA}
\lambda A(t) \leq A(\lambda t) \quad \hbox{for $t \geq 0$.}
\end{equation}
As a consequence, if $\lambda \geq 1$, then
\begin{equation}\label{lambdaA-1}
 A^{-1}(\lambda r) \leq \lambda A^{-1}(r) \quad \hbox{for $r \geq 0$.}
\end{equation}
\par\noindent
 A Young function $A$ is said to satisfy the
$\Delta_2$-condition
 if there exists a positive constant $C$ such that
\begin{align}\label{delta2}
A(2t)\leq CA(t)\quad \textrm{for \,\,} t\geq 0.
\end{align}
 We say that $A$ satisfies the $\nabla_2$-condition if
there exists a constant $C>2$ such that
\begin{equation}\label{nabla2}
A(2t) \geq C A(t) \quad \textrm{for \,\,} t\geq 0.
\end{equation}
 If $A$ is finite-valued and \eqref{delta2} just
holds for $t \geq t_0$ for some $t_0
>0$, then $A$ is said to satisfy the  $\Delta_2$-condition near
infinity. Similarly, if \eqref{nabla2} holds for $t \geq t_0$ for
some $t_0
>0$, then $A$ is said to satisfy the $\nabla_2$-condition near
infinity. We shall also write $A\in \Delta _2$ [$A \in \nabla _2$]
to denote that $A$ satisfies the  $\Delta_2$-condition
[$\nabla_2$-condition].
\par\noindent
One has that $A \in \Delta_2$ [near infinity] if and only if
$\widetilde A \in \nabla_2$ [near infinity].
\par\noindent
A Young function $A$ is said to dominate another Young function $B$
[near infinity] if there exists a positive constant  $C$
\begin{equation}\label{B.5bis}
B(t)\leq A(C t) \qquad \textrm{for \,\,\,$t\geq 0$\,\, [$t\geq t_0$
\,\, for some $t_0>0$]\,.}
\end{equation}
The functions $A$ and $B$ are called equivalent [near infinity] if
they dominate each other [near infinity]. We shall write $A \approx
B$ to denote such equivalence.
\par
Let $\Omega$ be a measurable subset of $\R^n$, and let $A$ be a
Young function. The Luxemburg norm associated with  $A$  is defined
as
\begin{align*}
\|u\|_{L^A(\Omega)}=\inf\left\{\lambda>0:\,\,\int_\Omega
A\Big(\frac{|u(x)|}{\lambda}\Big)\,dx\leq 1\right\}
\end{align*}
for any measurable function $u: \Omega \to \mathbb R$. The
collection of all  functions $u$ for which such norm is finite is
called the  Orlicz space $L^A(\Omega)$, and  is a Banach function
space. \\ A H\"older type inequality in Orlicz spaces takes the form
\begin{equation}\label{B.5}
\|v\|_{L^{\widetilde{A}}(\Omega)}\leq \sup _{u \in L^A(\Omega)}
\frac {\int _\Omega u(x) v(x)\, dx}{\|u\|_{L^A(\Omega)}} \leq 2
\|v\|_{L^{\widetilde{A}}(\Omega)}
\end{equation}
for every $v\in L^{\widetilde{A}}(\Omega)$.
 \\  Assume that
 $|\Omega|< \infty$ [$|\Omega|=\infty$], where $\|\, \cdot \,\|$ denotes Lebesgue
 measure, and let $A$ and $B$ be Young functions. Then
\begin{equation}\label{B.6}
L^A(\Omega)\to L^B(\Omega),
\end{equation}
if and only if $A$ dominates $B$ near infinity [globally]. The norm
of the embedding \eqref{B.6} depends on the
 constant $C$ appearing in
\eqref{B.5bis} if $A$ dominates $B$ globally. When $|\Omega | <
\infty$, and $A$ dominates $B$ just near infinity,  the embedding
constant also depends on $A$, $B$, $t_0$ and $|\Omega|$.
\\ The decreasing rearrangement $u^{\ast}: [0, \infty ) \to [0,
\infty]$ of a measurable function $u: \Omega \to \mathbb R$ is the
(unique) non-increasing, right-continuous function  which is
equimeasurable with $u$. Thus,
\begin{equation*}
u^{\ast}(s) = \inf \{t\geq 0: |\{x\in \Omega: |u(x)|>t \}|\leq s \} \qquad
\mathrm{for} \, s\geq 0.
\end{equation*}
The equimeasurability of $u$ and $u^*$ implies that
\begin{equation}\label{B.3}
\|u \|_{L^A(\Omega)} = \|u^{\ast} \|_{L^A(0,|\Omega|)}
\end{equation}
for every $u\in L^A(\Omega)$.
\par
The Lebesgue spaces $L^p(\Omega)$, corresponding to the choice
$A(t)=t^p$, if $p \in [1, \infty)$, and $A(t)= \infty \chi_{(1,
\infty)}(t)$, if $p=\infty$, are a basic example of Orlicz spaces.
Other customary instances of Orlicz spaces are provided by the
Zygmund spaces $L^p\log^\alpha L(\Omega)$, where either $p>1$ and
$\alpha \in \R$, or $p=1$ and $ \alpha \geq 0$, and by the
exponential spaces $\exp L^\beta (\Omega)$, where $\beta >0$. Here,
and in what follows, the notation $A (L)(\Omega)$ stands for the
Orlicz space associated with a Young function equivalent to the
function $A$ near infinity. Similar notations will be employed for
other function spaces, built upon Young functions, to be defined
below.
\par\noindent The Orlicz space $L^{A}(\Omega ,\R ^n)$
 of $\R^n$-valued measurable functions on $\Omega$ is defined
as $L^{A}(\Omega , \R^n)  = (L^{A}(\Omega))^n$, and is equipped with
the norm given by $\|\bfu \|_{L^{A}(\Omega , \R^n)}=\|\,|\bfu |\,
\|_{L^A(\Omega)}$ for $\bfu \in L^{A}(\Omega , \R^n)$. The  Orlicz
space $L^{A}(\Omega , \mathbb R^{n\times n})$ of  matrix-valued
measurable functions on $\Omega$ is defined analogously.

\medskip
\par

Assume now that $\Omega$ is an open set.  The Orlicz-Sobolev space
$W^{1,A}(\Omega)$  is the set of all weakly differentiable functions
in $L^A(\Omega)$
 whose  gradient belongs
to $L^A(\Omega,\rn)$. The alternate notation $W^1L^A(\Omega)$ for
$W^{1,A}(\Omega)$ will also be used when convenient. The space
$W^{1,A}(\Omega)$  is a Banach space endowed with the norm
\begin{align*}
\|u\|_{W^{1,A}(\Omega)}=\|u\|_{L^A(\Omega)}+\|\nabla
u\|_{L^A(\Omega, \rn)}.
\end{align*}
We also define
 \begin{align*} W^{1,A}_0(\Omega) = \{u \in
W^{1,A}(\Omega) :\,\, & \hbox{the continuation of $u$ by $0$ outside
$\Omega$} \\  & \qquad \quad \hbox{is weakly differentiable in $\R
^n$}\}.
\end{align*}
In the case when $A(t)=t^p$ for some $p \geq 1$, and $\partial
\Omega$ is regular enough, such definition of $W^{1,A}_0(\Omega)$
can be shown to reproduce the usual space $W^{1,p}_0(\Omega)$
defined as the closure in $W^{1,p}(\Omega )$ of the space $C^\infty
_0(\Omega )$ of smooth compactly supported functions in $\Omega$. In
general, the set of smooth bounded functions is dense in $L^A(\Omega
)$ only if $A$ satisfies the $\Delta _2$-condition (just near
infinity when $|\Omega|< \infty$). Thus, for arbitrary $A$, our
definition of $W^{1,A}_0(\Omega)$ yields a space which can be larger
than the closure of $C^\infty _0(\Omega )$ in $W^{1,A}_0(\Omega)$
even for a set $\Omega$ with a smooth boundary. On the other hand,
if $\Omega$ is a bounded Lipschitz domain,
 then $W^{1,A}_0(\Omega)= W^{1,A}(\Omega)\cap W^{1,1}_0(\Omega),$
where $W^{1,1}_0(\Omega)$ is defined as usual. Recall that an open
set $\Omega$ is called a Lipschitz domain if  there exists a
neighborhood $\mathcal U$  of each point of $\partial \Omega$ such
that $\Omega \cap \mathcal U$ is  the subgraph of a Lipschitz
continuous function of $n-1$ variables. An open set $\Omega$ is said
to have the cone property if there exists a finite cone $\Lambda$
such that each point of $\Omega$ is the vertex of a finite cone
contained in $\Omega$ and congruent to $\Lambda$. Moreover, an open
set $\Omega$ is said to be starshaped with respect to a ball
$\mathcal B \subset \Omega$ if it is starshaped with respect to
every point in $\mathcal B$. Clearly, any bounded open set which is
starshaped with respect to a ball is a Lipschitz domain, and any
bounded Lipschitz domain has the cone property.
\par
The Orlicz-Sobolev space $W^{1, A} (\Omega , \R^n)$ of $\R^n$-valued
functions is defined as $W^{1, A}(\Omega , \R^n) = \big(W^{1,
A}(\Omega )\big)^n$, and equipped with the norm $\|\bfu\|_{W^{1, A}
(\Omega ,\R^n)} =  \|\bfu \|_{L^A(\Omega ,\rn)} + \|\nabla
\bfu\|_{L^A(\Omega , \rmat)}$. The space $W^{1, A}_0(\Omega , \R^n)$
is defined accordingly.

\smallskip
\par

We next denote by $E^A(\Omega, \rn)$, or by   $EL^A(\Omega, \rn)$,
the space of those functions $\bfu \in L^A(\Omega , \rn)$ whose
distributional symmetric gradient
$$\mathcal E \bfu = \tfrac 12(\nabla \bf u + (\nabla \bf u)^{T})$$
belongs to $L^A(\Omega , \rmat)$. Here, $\lq\lq (\, \cdot \,)^T"$
stands for transpose.  $E^A(\Omega, \rn)$ is a Banach space equipped
with the norm
\begin{equation}\label{normE}
\|\bfu\|_{E^{A} (\Omega ,\R^n)} =  \|\bfu \|_{L^A(\Omega ,\rn)} +
\|\mathcal E \bfu\|_{L^A(\Omega , \rmat)}.
\end{equation}
The subspace $E^A_0(\Omega, \rn)$ is defined as the set of those
functions in $E^A(\Omega, \rn)$ whose continuation by $0$ outside
$\Omega$ belongs to $E^A(\rn, \rn)$.
\\
The kernel of the operator $\mathcal E$, in any connected open set
$\Omega$ in $\rn$,  is known to agree with the space
$$\mathcal R = \{{\bf v} :\rn \to \rn: {\bf v}(x)= \bfb+ {\bf Q}x \,\,
\hbox{for some $\bfb \in \rn$ and  ${\bf Q}\in \rmat$ such that
${\bf Q}=-{\bf Q}^T$}\},$$ see e.g. \cite[Lemma 1.1, Chapter 1]{Te}.

\smallskip
\par

The notation $E^{D,A}(\Omega , \rn)$ is devoted to the space of
those functions $\bfu \in L^A(\Omega , \rn)$ whose trace-free
distributional symmetric gradient
$$\mathcal E^D  \bfu = \mathcal E \bfu - \tfrac {{\rm tr} (\mathcal E \bfu)} n I $$
belongs to $L^A(\Omega , \rmat)$.  Here, $I$ denotes the identity
matrix, and ${\rm tr} (\mathcal E \bfu)$ the trace of the matrix
$\mathcal E \bfu$. The space $E^{D,A}(\Omega , \rn)$, which will
also be occasionally denoted by $E^{D}L^A(\Omega, \rn)$, is also a
Banach space equipped with the norm
\begin{equation}\label{normED}
\|\bfu\|_{E^{D,A} (\Omega ,\R^n)} =  \|\bfu \|_{L^A(\Omega ,\rn)} +
\|\mathcal E^D \bfu\|_{L^A(\Omega , \rmat)}.
\end{equation}
The definition of the subspace~$E^{D,A}_0(\Omega, \rn)$ of
  $E^{D,A}(\Omega, \rn)$ parallels those of $W^{1,A}_0(\Omega, \rn)$ and $E^A(\Omega , \rn)$.
\par\noindent
The kernel of the operator $\mathcal E^D$, in any  connected open
set $\Omega$ in $\rn$, $n \geq 3$, is  the direct sum
$$\Sigma  = \mathcal D\oplus\mathcal R\oplus\mathcal S,$$
where
\begin{align*}
\mathcal D&=\set{\bfv:\rn\rightarrow\rn:\,\,\bfv(x)=\rho x\,\,\text{for some}\,\,\rho\in\R},\\
\mathcal S&=\set{\bfv:\rn\rightarrow\rn:\,\,\bfv(x)=2(\bfa\cdot
x)x-|x|^2\bfa\,\,\text{for some}\,\,\bfa\in\R^n},
\end{align*}
see e.g. \cite[Prop. 2.5]{Sc2}.

\section{Main results}\label{main}

Our characterization of the Young functions $A$ and $B$ supporting
trace-free Korn-type inequalities between the Orlicz spaces $L^A$
and $L^B$ amounts to the balance conditions:
\begin{subequations}
  \begin{align}\label{1.1}
    t\int _{t_0} ^t \frac{ B(s)}{s^2} \,ds \leq A(ct) \qquad
    \hbox{for $t \geq t_0$,}
  \end{align}
  and
  \begin{align}\label{1.2}
    t \int _{t_0}^t \frac{ \widetilde A(s)}{s^2} \,ds \leq \widetilde
    B(ct) \qquad \hbox{for $t \geq t_0$,}
  \end{align}
\end{subequations}
for some  constants $c>0$ and $t_0\geq0$.

\smallskip
\par\noindent
The result for functions vanishing on the boundary of their domain
reads as follows.

\begin{theorem}\label{thm:main0} {\bf [Trace-free Korn inequalities
in $E_0^{D,A}(\Omega,\rn)$]} Let $\Omega$ be an open bounded set in
$\rn$, $n\geq3$.
 Let $A$
and $B$ be Young functions. The following facts are equivalent.
\par\noindent (i) Inequalities \eqref{1.1} and \eqref{1.2} hold.
\par\noindent (ii)   $E_0^{D,A}(\Omega,\rn)\subset W_0^{1,B}(\Omega,\rn)$, and there exists a constant $C$ such
that
\begin{equation}\label{1.3.0}
  \| \nabla \bfu\|_{L^B(\Omega , \rmat)} \leq C \| \mathcal E^D \bfu\|_{L^A(\Omega
, \rmat)}
\end{equation}
 for every $\bfu\in E^{D,A}_0(\Omega , \rn)$.
\par\noindent (iii)  $E_0^{D,A}(\Omega,\rn)\subset W_0^{1,B}(\Omega,\rn)$, and there exist constants $C$ and $C_1$ such
that
\begin{equation}\label{1.3.1}
\int_\Omega B(|\nabla\bfu|)\,dx\leq\,C_1+\int_\Omega A(C|\mathcal E^D\bfu|)\,dx
\end{equation}
 for every $\bfu\in E^{D,A}_0(\Omega , \rn)$.
\end{theorem}

\begin{remark}\label{remarkrn}
{\rm A close inspection of the proof of Theorem \ref{thm:main0}
reveals that inequality \eqref{1.3.1} holds with $C_1=0$ if and only
if conditions \eqref{1.1} and \eqref{1.2} are fulfilled  with
$t_0=0$. When $\Omega = \rn$,  these conditions with $t_0=0$ turn
out to be equivalent to inequalities \eqref{1.3.0} and
\eqref{1.3.1}, with $C_1=0$. In fact,  if \eqref{1.3.1} holds with
$\Omega = \rn$ for some $C_1$, then it also holds with $C_1=0$. This
follows from  a scaling argument, based on replacing any trial
function $\bfu(x)$ in \eqref{1.3.1} with $R\bfu(x/R)$ for $R>0$, and then
letting $R\to \infty$.}
\end{remark}

\par\noindent
Inequalities without  boundary conditions are the object of the next
theorem.

\begin{theorem}\label{thm:main} {\bf [Trace-free Korn inequalities
in $E^{D,A}(\Omega,\rn)$]} Let $\Omega$ be a   bounded  connected
open set with the cone property in $\rn$, $n\geq3$.
 Let $A$
and $B$ be Young functions. The following facts are equivalent.
\par\noindent (i) Inequalities \eqref{1.1} and \eqref{1.2} hold.
\par\noindent (ii)  $E^{D,A}(\Omega,\rn)\subset W^{1,B}(\Omega,\rn)$, and there exists a constant $C$ such that
\begin{equation}\label{1.3main}
\inf _{\bfw \in \Sigma}  \| \nabla \bfu - \nabla \bfw\|_{L^B(\Omega
, \rmat)} \leq C \| \mathcal E^D \bfu\|_{L^A(\Omega , \rmat)}
\end{equation}
 for every $\bfu\in E^{D,A}(\Omega , \rn)$.
\par\noindent (iii)  $E^{D,A}(\Omega,\rn)\subset W^{1,B}(\Omega,\rn)$, and there exist constants $C$ and $C_1$ such
that
\begin{equation}\label{1.3.2}
\inf_{\bfw\in\Sigma}\int_\Omega B(|\nabla\bfu-\nabla\bfw|)\,dx\leq\,C_1+\int_\Omega A(C|\mathcal E^D\bfu|)\,dx
\end{equation}
 for every $\bfu\in E^{D,A}(\Omega , \rn)$.
\end{theorem}

\begin{remark}\label{remark1march}
{\rm Similarly to \eqref{1.3.1}, if conditions
\eqref{1.1}--\eqref{1.2} are fulfilled with $t_0=0$, then inequality
\eqref{1.3.2} holds with $C_1=0$.}
\end{remark}

\begin{remark}\label{remark2march}
{\rm If  either \eqref{1.1} or \eqref{1.2} is in force,  then $A$
dominates $B$ near infinity, or globally, according to whether
$t_0>0$ or $t_0=0$ \cite[Proposition 3.5]{Ci2}. Moreover, inequality
\eqref{1.1} holds with $B=A$ for some $t_0>0$ [resp. for $t_0=0$],
if and only if $A \in \nabla _2$ near infinity [resp. globally], and
inequality \eqref{1.2} holds with $B=A$ for some $t_0>0$ [for
$t_0=0$] if and only if $A \in \Delta _2$ near infinity [globally]
\cite[Theorem 1.2.1]{KK}.
\\ Thus, Theorems \ref{thm:main0} and \eqref{thm:main} recover the
fact that inequalities \eqref{1.3.0}--\eqref{1.3.1} and
\eqref{1.3main}--\eqref{1.3.2} hold with
 $B=A$ if and only if $A \in \Delta_2 \cap \nabla _2$ near infinity.
}
\end{remark}

Hereafter, we present some inequalities for functions in  spaces
$E^{D,A}_0(\Omega , \rn)$ of logarithmic or exponential type, which
follow from
 Theorem \ref{thm:main0} and Remark
\ref{remark2march}.  Analogues for $E^{D,A}(\Omega , \rn)$ hold
owing  to Theorem \ref{thm:main}, provided that $\Omega$ fulfils the
assumptions of the latter. In the following examples $\Omega$
denotes a bounded open set in $\rn$, $n \geq 3$.

\begin{example}\label{ex1}
{\rm   If $p>1$ and $\alpha \in \R$, then
\begin{equation}\label{ex1.1}
\| \nabla \bfu \|_{L^p(\log  L)^\alpha (\Omega , \rmat)} \leq C \|
\mathcal E^D \bfu\|_{L^p(\log  L)^\alpha (\Omega , \rmat)}
\end{equation}
for every $\bfu \in E_0^{D}L^p(\log  L)^\alpha (\Omega , \rn)$. If
$\alpha \geq 0$, then
\begin{equation}\label{ex1.2}
\| \nabla \bfu \|_{L (\log  L)^\alpha (\Omega , \rmat)} \leq C \|
\mathcal E^D \bfu\|_{L (\log  L)^{\alpha +1} (\Omega , \rmat)}
\end{equation}
for every $\bfu \in E_0^{D}L (\log  L)^{\alpha +1} (\Omega , \rn)$.}
\end{example}

\medskip

\begin{example}\label{exrn}
{\rm  Let  $B$ be a Young functions such that
$$B(t) \approx \begin{cases} t^q \Big(\log \frac 1t\Big)^{-\beta} &
\hbox{near $0$} \\
t^p \Big(\log t\Big)^\alpha & \hbox{near $\infty$,}\end{cases}$$
where either $q > 1$ and $\beta \in \R$, or $q=1$ and $\beta
>1$, and either $p > 1$ and $\alpha \in \R$, or  $p=1$ and $\alpha \geq
0$. Assume that $A$ is another Young function fulfilling
$$A(t) \approx \begin{cases} B(t) & \hbox{if $q>1$} \\
t \Big(\log \frac 1t\Big)^{1-\beta} & \hbox{if $q=1$,}\end{cases}
$$
near $0$, and
$$A(t) \approx \begin{cases} B(t) & \hbox{if $p>1$} \\
t \Big(\log t\Big)^{1+\alpha} & \hbox{if $p=1$,}\end{cases}
$$
near  infinity. Then
\begin{equation}\label{exrn.1}
\| \nabla \bfu \|_{L^B (\rn , \rmat)} \leq C \| \mathcal E^D
\bfu\|_{L^A (\rn , \rmat )}
\end{equation}
for every compactly supported function $\bfu \in E^{D,A}(\rn ,
\rn)$, with $n \geq 3$.}
\end{example}

\medskip

\begin{example}\label{ex2}
{\rm Assume that $p>1$ and $\alpha \in \R$. Then
\begin{equation}\label{ex2.1}
\| \nabla \bfu \|_{L^p(\log  \log  L)^\alpha (\Omega , \rmat)} \leq
C \| \mathcal E^D \bfu\|_{L^p(\log  \log  L)^\alpha (\Omega ,
\rmat)}
\end{equation}
for every $\bfu \in E_0^{D}L^p(\log  \log  L)^\alpha (\Omega ,
\rn)$. If $\alpha \geq 0$, then
\begin{equation}\label{ex2.2}
\| \nabla \bfu \|_{L (\log  \log  L)^\alpha (\Omega , \rmat)} \leq C
\| \mathcal E^D \bfu\|_{L \log  L (\log  \log  L)^{\alpha} (\Omega ,
\rmat)}
\end{equation}
for every $\bfu \in E_0^{D}L \log  L (\log  \log  L)^{\alpha}
(\Omega , \rn)$.}
\end{example}

\medskip

\begin{example}\label{ex3}
{\rm Assume that $\beta
>0$. Then
\begin{equation}\label{ex3.1}
\| \nabla \bfu \|_{\exp  L^{\frac{\beta}{\beta +1}} (\Omega ,
\rmat)} \leq C \| \mathcal E^D \bfu\|_{\exp  L^{\beta} (\Omega ,
\rmat)}
\end{equation}
for every $\bfu \in E_0^{D}\exp  L^{\beta} (\Omega , \rn)$.}
\end{example}

\medskip

\begin{example}\label{ex5}
{\rm One has that
\begin{equation}\label{Linf}
\| \nabla \bfu \|_{\exp  L (\Omega , \rmat)} \leq C \| \mathcal E^D
\bfu\|_{  L^{\infty} (\Omega , \rmat)}
\end{equation}
for every $\bfu \in E_0^{D}   L^{\infty} (\Omega , \rn)$.}
\end{example}

\medskip

\begin{example}\label{ex4}
{\rm Assume that $a>0$ and $\beta >1$. Then
\begin{equation}\label{ex4.1}
\| \nabla \bfu \|_{\exp  (a (\log  L)^\beta )  (\Omega , \rmat)}
\leq C \| \mathcal E^D \bfu\|_{\exp  \big(a \big(\log
\frac{L}{(\log L)^{\beta -1}}\big)^\beta \big) (\Omega , \rmat)}
\end{equation}
for every $\bfu \in E_0^{D}\exp  \Big(a \Big(\log  \frac{L}{(\log
L)^{\beta -1}}\Big)^\beta \Big) (\Omega , \rn)$.}
\end{example}

\medskip

The  necessity of conditions \eqref{1.1} and \eqref{1.2} in our
results about trace-free Korn inequalities goes through a proof of
 their necessity
in the Orlicz-Korn inequality for the plain symmetric gradient. The
sufficiency of \eqref{1.1} and \eqref{1.2} for the latter inequality
was established in \cite{Ci2}. A comprehensive statement,
summarizing necessary and sufficient conditions for the Korn
inequality in Orlicz spaces, reads as follows.

\begin{theorem}\label{thm:plainkorn}  {\bf [Korn inequalities in $E^{A}_0(\Omega , \rn)$ and $E^{A}(\Omega ,
\rn)$]}\, {\rm (see also \cite[Theorems 3.1 and 3.3]{Ci2})}
 Let $A$
and $B$ be Young functions. The following facts are equivalent.
\par\noindent (i) Inequalities \eqref{1.1} and \eqref{1.2} hold.
\par\noindent (ii)  Given an bounded open set $\Omega$ in $\rn$, $n \geq 2$, the inclusion $E_0^{A}(\Omega,\rn)\subset W_0^{1,B}(\Omega,\rn)$ holds,
and there exists a constant $C$ such
that
\begin{equation}\label{1.3.plainkorn}
  \| \nabla \bfu\|_{L^B(\Omega , \rmat)} \leq C \| \mathcal E \bfu\|_{L^A(\Omega
, \rmat)}
\end{equation}
 for every $\bfu\in E^{A}_0(\Omega , \rn)$.
\par\noindent (iii)  Given a
bounded connected open set with the cone property in $\rn$, $n\geq
2$, the inclusion $E^{A}(\Omega,\rn)\subset W^{1,B}(\Omega,\rn)$
holds, and there exists a constant $C$ such that
\begin{equation}\label{1.3plainmainkorn}
\inf _{\bfv \in \mathcal R}  \| \nabla \bfu - \nabla
\bfv\|_{L^B(\Omega , \rmat)} \leq C \| \mathcal E \bfu\|_{L^A(\Omega
, \rmat)}
\end{equation}
 for every $\bfu\in E^{A}(\Omega , \rn)$.
\end{theorem}

\medskip
\par\noindent
Conditions  \eqref{1.1} and \eqref{1.2} also appear in an inequality
for negative Orlicz-Sobolev norms recently established in
\cite{BrCi}. Let $A$ be a Young function. The negative
Orlicz-Sobolev norm of the distributional gradient of a function $u
\in L^1(\Omega)$ can be defined as
\begin{align}\label{negorlicz}
\|\nabla u\|_{W^{-1,A}(\Omega, \rn)} =\sup_{\bfvarphi\in
C^\infty_0(\Omega, \R^n)} \frac{\int_\Omega u\,{\rm div}\, \bfvarphi
\,dx}{\|\nabla \bfvarphi\|_{L^{\widetilde{A}}(\Omega, \R^{n\times
n})}}.
\end{align}
This definition, introduced in \cite{BrCi}, is an Orlicz space
version of negative norms for classical Sobolev spaces which
goes back to Ne\v{c}as \cite{Ne}. He showed that, if $\Omega$ is
regular enough,  and $1 < p < \infty$, then the $L^p(\Omega)$ norm
of any function with zero mean-value over $\Omega$ is equivalent to
the $W^{-1,p}(\Omega, \rn)$  norm of its gradient, defined as in
\eqref{negorlicz} with $L^{\widetilde{A}}(\Omega, \R^{n\times n}) =
L^{p'}(\Omega , \R^{n\times n})$, and $p'=\tfrac p{p-1}$. Namely,
there exist positive constants $C_1$ and $C_2$ such that
$$
C_1\|u-u_\Omega\|_{L^p(\Omega)}\leq \|\nabla u\|_{W^{-1,p}(\Omega,
\rn)}\leq C_2\|u-u_\Omega\|_{L^p(\Omega)}
$$
for every $u\in L^1(\Omega)$, where $u_\Omega=\tfrac 1{|\Omega|}\int
_\Omega u\,dx$ the mean value of $u$ over $\Omega$.
\\ The inequality:
\begin{equation}\label{negtrivial}
\|\nabla u\|_{W^{-1,A}(\Omega,\rn)} \leq C
\|u-u_\Omega\|_{L^A(\Omega)}
\end{equation}
holds for every Young function $A$, for some absolute constant $C$,
and for every $u \in L^1(\Omega)$ \cite[Theorem 3.1]{BrCi}. Although
a reverse inequality fails in general, it can be restored provided
that the norm of $u-u_\Omega$ in $L^A(\Omega)$ is replaced with the
norm in some Orlicz space $L^B(\Omega)$, with $B$ fulfilling
\eqref{1.1} and \eqref{1.2}. This is also established in
\cite[Theorem 3.1]{BrCi}. The necessity of conditions \eqref{1.1}
and \eqref{1.2} for the relevant reverse inequality follows from
their necessity in Theorem \ref{thm:plainkorn}. Altogether, the
following result holds.

\begin{theorem}\label{negativenecsuf} {\bf [Negative norm inequalities]}\, {\rm (see also \cite[Theorem
3.1]{BrCi})}
 Let $A$ and $B$ be
Young functions. Let $\Omega$ be a  bounded connected open set with
the cone property in $\rn$, $n\geq 2$. There exists
  a constant $C$ such that
\begin{equation}\label{march1}
\|u-u_\Omega\|_{L^B(\Omega)}\leq C\|\nabla
u\|_{W^{-1,A}(\Omega,\rn)}
\end{equation}
for every $u \in L^1(\Omega)$ if and only if $A$ and $B$ satisfy
conditions \eqref{1.1} and \eqref{1.2}.
\end{theorem}

\bigskip
\par
The proof of   inequality  \eqref{march1} relies upon boundedness
properties of the gradient of the \Bogovskii \, operator. Given a
bounded open set $\Omega$, which is starshaped with respect to some
ball, and any smooth, nonnegative function $\omega$, compactly
supported in such ball and with integral equal to $1$, the
\Bogovskii \, operator $\mathcal B_\Omega$ is defined, according to
\cite{Bog}, as
\begin{align}\label{bog}
\mathcal B_\Omega f(x)=\int_{\Omega}
f(y)\bigg(\frac{x-y}{|x-y|^{n}}\int_{|x-y|}^\infty\omega\Big(y+ r
\frac{x-y}{|x-y|}\Big)\zeta^{n-1}\,dr \bigg)\,dy \quad \hbox{for $x
\in \Omega$,}
\end{align}
for every function $f \in C^{\infty}_{0, \bot}(\Omega)$. Here,
$C^{\infty}_{0, \bot}(\Omega)$ denotes the subspace of
$C^{\infty}_{0}(\Omega)$ of those functions with vanishing
mean-value
 on $\Omega$.  This operator  is customarily used to construct a solution
 to the divergence equation, coupled  with zero boundary conditions, inasmuch as
$\Div\mathcal B_\Omega f=f$.
\\ The boundedness of the   operator $\nabla \mathcal B_\Omega$
between Orlicz spaces $L^A(\Omega)$ and $L^B(\Omega,\rmat)$, under
assumptions \eqref{1.1} and \eqref{1.2}, is proved in
\cite[inequality (3.88)]{BrCi}. The necessity part of Theorem
\ref{negativenecsuf} allows to show that these assumptions are, in
fact, also necessary. In conclusion,  the following full
characterization holds.

\bigskip

\begin{theorem}\label{corbogov} {\bf [Boundedness properties of $\nabla \mathcal B_\Omega$]} \, {\rm [see also \cite[Theorem
3.6]{BrCi}]} Let $A$ and $B$ be Young functions. Let  $\Omega$ be
a bounded open set in $\rn$, $n \geq 2$, which is starshaped with
respect to a ball.  There exists a constant $C$ such that
\begin{equation}\label{bogov1}
\|\nabla \mathcal B_\Omega f\|_{L^B(\Omega , \rmat)} \leq C
\|f\|_{L^A(\Omega)}
\end{equation}
 for every $f \in C^{\infty}_{0,
\bot}(\Omega)$ if and only if $A$ and $B$ satisfy conditions
\eqref{1.1} and \eqref{1.2}.
\end{theorem}

\section{Representation formulas and trace-free Korn inequalities in Orlicz spaces}
\label{sec:Korn2suff}

We are concerned here with a proof of inequalities \eqref{1.3.0} and
\eqref{1.3main} under conditions \eqref{1.1} and \eqref{1.2}. The
former inequality, which involves functions vanishing on the
boundary of their domain, is the object of the first result.

\begin{theorem}\label{sufficiency0}  Let $\Omega $ be a  bounded open set in $\rn$, $n
\geq 3$. Let $A$ and $B$ be Young functions fulfilling conditions
\eqref{1.1} and \eqref{1.2}. Then
 $E^{D,A}_0(\Omega , \rn)
\subset W^{1,B}_0(\Omega , \rn)$, and inequality \eqref{1.3.0}
holds.
\end{theorem}

\smallskip

 The relevant inequality for arbitrary functions is established in
 the next theorem.

\smallskip

\begin{theorem}\label{sufficiency}
Let $\Omega $ be a  bounded connected open set with the cone
property in $\rn$, $n \geq 3$. Assume that  $A$ and $B$ are Young
functions fulfilling conditions \eqref{1.1} and \eqref{1.2}. Then
 $E^{D,A}(\Omega , \rn)
\subset W^{1,B}(\Omega , \rn)$,  and inequality \eqref{1.3main}
holds.
%
%
\end{theorem}

The proofs of Theorem~\ref{sufficiency0} and
Theorem~\ref{sufficiency} are split into several lemmas, and are
accomplished at the end of this section. We begin with  Lemma \ref{repr}, whose objective is to show that   the full gradient can be represented as a singular integral of $\mathcal E ^D$, plus  some weaker  terms, also depending on $\mathcal E ^D$. In Lemma \ref{febbraio100}, a pointwise estimate, in rearrangement form, is established for the relevant singular integral operator. This  reduces   the question of the validity of  a trace-free Korn-type inequality in Orlicz spaces to that of a considerably simpler one-dimensional Hardy inequality in the same spaces. General criteria for the Hardy inequalities that come into play are stated in Lemma \ref{hardy}.

\medskip
\par\noindent

\begin{lemma}\label{repr}
Let $\Omega$ be a bounded open set in $\rn$, $n \geq 3$, which is
starshaped with respect to a ball. Let $A(t)$ be a Young function
 which dominates the function $t\log(1+t)$ near infinity. Assume that $\bfu
\in E^{D,A}(\Omega, \rn)$. Then $\bfu \in W^{1,1}_{\rm loc}(\Omega ,
\rn)$, and
\begin{equation}\label{july1}
\frac{\partial u_h}{\partial x_k}(x) = P_{hk} (x) + \sum _{i,j=1}^n
\int _{\Omega} \mathcal E ^D_{ij}(\bfu)(y) K_{ijhk}(x, y)\, dy +
 \sum _{i,j=1}^n C_{ijhk}\mathcal E ^D_{ij}(\bfu)(x)\quad \hbox{for a.e. $x \in
 \Omega$,}
 \end{equation}
where  $u_h$ denotes the $h$-th component of $\bfu$, $\mathcal E
^D_{ij}(\bfu)$   the $ij$ entry of the matrix $\mathcal E ^D(\bfu)$,
$P_{hk}$ are polynomials of degree one, $C_{ijhk}$ are constants,
and $K_{ijhk}:\Omega\times\rn\rightarrow\mathbb R$ are kernels of
the form $K_{ijhk}(x,y) = N(x, y-x)$ for some function
$N:\Omega\times (\rn  \setminus \{0\}) \rightarrow\mathbb R$,
depending on $\Omega$ and on $i,j,j,k$, and enjoying the following
properties:
\begin{equation}\label{N1}
N(x,\lambda z)= \lambda ^{-n} N(x, z) \quad \hbox{for $x \in
\Omega$, $z \in \rn \setminus \{0\}$, $\lambda >0$;}
\end{equation}
\begin{equation}\label{N2}
\int _{\mathbb S^{n-1}}N(x, z)\, d\mathcal H^{n-1} (z) =0 \quad
\hbox{for  $x \in \Omega$;}
\end{equation}
where $\mathcal H^{n-1}$ denotes the surface measure on $\mathbb S^{n-1}$;\\
 for every $p \in [1,
\infty)$ there exists a constant $C$ such that
\begin{equation}\label{N3}
\int _{\mathbb S^{n-1}}|N(x, z)|^p\, d\mathcal H^{n-1} (z) \leq C
\quad \hbox{for  $x \in \Omega$;}
\end{equation}
there exists a constant $C$ such that
\begin{equation}\label{N4}
|N(x, y-x)| \leq \frac C{|x-y|^{n}} \quad \hbox{for $x\neq y$,}
\end{equation}
\begin{equation}\label{N5}
|N(x, y-x) - N(z, y-z)| \leq C \frac{|z-x|}{|y-x|^{n+1}}\quad \hbox{for $x\neq y$ and $2|x-z| < |x-y|$,}
\end{equation}
\begin{equation}\label{N6}
|N(y, x-y) - N(y, z-y)| \leq C \frac{|z-x|}{|y-x|^{n+1}}\quad
\hbox{for $x\neq y$ and $2|x-z| < |x-y|$}.
\end{equation}
 \end{lemma}
%
\par\noindent
{\bf Proof}. The representation formula \cite[Equation (2.43)]{Resh}
tells us that, if  $\bfu \in C^\infty (\Omega , \rn)$, then
\begin{equation}\label{repres}
\bfu (x) = P\bfu (x) + R (\mathcal E ^D \bfu) (x) \quad \hbox{for $x
\in \Omega$}.
\end{equation}
Here, for each $i=1, \dots , n$,
\begin{equation}\label{repres1}
(P\bfu )_i (x) = \sum _{0 \leq |\alpha | \leq 2} x^\alpha \sum
_{k=1}^n \int _\Omega u_k (y) H_{ik\alpha}(y) \, dy \quad \hbox{for
$x \in \Omega$},
\end{equation}
where $u_k$ denotes the $k$-th component of $\bfu$, $(P\bfu )_i$ the
$i$-th component of $P\bfu$, and the functions $H_{ik\alpha}\in
C^\infty _0(\Omega)$ are such that $P\bfu \in \Sigma $ for every
$\bfu \in C^\infty (\Omega , \rn)$. The expression $x^\alpha$
denotes a polynomial of the form $x_1^{\alpha_1}x_2^{\alpha
_2}\cdots x_n^{\alpha _n}$, where $\alpha = (\alpha _1, \dots ,
\alpha _n)$ is a multi-index of length  $|\alpha | = \alpha _1+
\cdots + \alpha _n$. Moreover,
\begin{equation}\label{repres2} R (\mathcal E ^D \bfu)_i (x) =
\sum _{k,j} \int _\Omega (\mathcal E ^D\bfu)_{kj} (y) R_{ikj}(x,
y)\, dy \quad \hbox{for $x \in \Omega$},
\end{equation}
where $R (\mathcal E ^D \bfu)_i$ is the $i$-th  component of $R
(\mathcal E ^D \bfu)$, and the kernels $R_{ikj}: \Omega\times\Omega
\setminus \{x=y\}\rightarrow\mathbb R$ are linear combinations, with
constant coefficients, of functions of the form
\begin{equation}\label{repres10}
(x_h - y_h) K(x,y),
\end{equation}
for some $h =1, \dots, n$, or
\begin{equation}\label{repres11}
\frac{\partial }{\partial {y_\ell}}\Big((x_h - y_h)(x_m - y_m)
K(x,y)\Big),
\end{equation}
for some $h, m, \ell =1, \dots, n$,
\begin{equation}\label{repres12}
\frac{\partial^2 }{\partial {y_\ell}\partial y_\kappa}\Big((x_h -
y_h)(x_m - y_m)(x_\iota - y_\iota) K(x,y)\Big),
\end{equation}
for some $h, m, \iota, \ell, \kappa =1, \dots, n$,
\begin{equation}\label{repres13}
K(x,y)= \frac{1}{|x-y|^n} \int _{|x-y|}^\infty \varphi \bigg(x +
\frac{y-x}{|y-x|}r\bigg) r^{n-1}\, dr\quad\text{for
}(x,y)\in\Omega\times\Omega, \, x \neq y,
\end{equation}
and  $\phi$ is any function in  $C^\infty _0 (\Omega)$. Note the
alternative formula:
\begin{equation}\label{repres14}
K(x,y)=  \int _{1}^\infty \varphi (x + (y-x)r) r^{n-1}\,
dr\quad\text{for }(x,y)\in\Omega\times\Omega, \, x \neq y.
\end{equation}
Making use of \eqref{repres14} in \eqref{repres11} and
\eqref{repres12}, and differentiating shows that the kernels $R_{kj}$
are linear combinations of functions of the form
\begin{equation}\label{repres15}
 \frac{1}{|x-y|^{n-1}} \frac {x_h-y_h}{|x-y|} \int _{|x-y|}^\infty \varphi \bigg(x +
\frac{y-x}{|y-x|}r\bigg) r^{n-1}\, dr,
\end{equation}
or
\begin{equation}\label{repres16}
 \frac{1}{|x-y|^{n-1}} \frac { x_h-y_h }{|x-y|}\frac{ x_m - y_m }{|x-y|} \int _{|x-y|}^\infty \frac{\partial \varphi}{\partial z_\ell} \bigg(x +
\frac{y-x}{|y-x|}r\bigg) r^{n}\, dr,
\end{equation}
or
\begin{equation}\label{repres17}
 \frac{1}{|x-y|^{n-1}} \frac {x_h-y_h }{|x-y|}\frac{ x_m - y_m }{|x-y|}\frac{ x_\iota - y_\iota }{|x-y|} \int _{|x-y|}^\infty
 \frac{\partial ^2 \varphi }{\partial z_\ell \partial z_\kappa} \bigg(x +
\frac{y-x}{|y-x|}r\bigg) r^{n+1}\, dr.
\end{equation}
In turn, these functions can be rewritten as
\begin{multline}\label{repres15bis}
 \frac{1}{|x-y|^{n-1}} \frac { x_h-y_h }{|x-y|} \int _{0}^\infty \varphi \bigg(x +
\frac{y-x}{|y-x|}r\bigg) r^{n-1}\, dr \\ - \frac{1}{|x-y|^{n-1}}
\frac { x_h-y_h }{|x-y|} \int _0^{|x-y|} \varphi \bigg(x +
\frac{y-x}{|y-x|}r\bigg) r^{n-1}\, dr,
\end{multline}
\begin{multline}\label{repres16bis}
 \frac{1}{|x-y|^{n-1}} \frac { x_h-y_h }{|x-y|}\frac{ x_m - y_m }{|x-y|} \int _{0}^\infty  \frac{\partial \varphi}{\partial z_\ell}\bigg(x +
\frac{y-x}{|y-x|}r\bigg) r^{n}\, dr \\
- \frac{1}{|x-y|^{n-1}} \frac { x_h-y_h }{|x-y|}\frac{ x_m - y_m
}{|x-y|} \int _0^{|x-y|}  \frac{\partial \varphi }{\partial z_\ell}
\bigg(x + \frac{y-x}{|y-x|}r\bigg) r^{n}\, dr ,
\end{multline}
\begin{multline}\label{repres17bis}
 \frac{1}{|x-y|^{n-1}} \frac { x_h-y_h }{|x-y|}\frac{ x_m - y_m }{|x-y|}\frac{ x_\iota - y_\iota }{|x-y|} \int _{0}^\infty
  \frac{\partial ^2 \varphi}{\partial z_\ell \partial z_\kappa} \bigg(x +
\frac{y-x}{|y-x|}r\bigg) r^{n+1}\, dr
\\
- \frac{1}{|x-y|^{n-1}} \frac { x_h-y_h }{|x-y|}\frac{ x_m - y_m
}{|x-y|}\frac{ x_\iota - y_\iota }{|x-y|} \int _0^{|x-y|}
 \frac{\partial ^2 \varphi}{\partial z_\ell \partial z_\kappa}\bigg(x + \frac{y-x}{|y-x|}r\bigg)
r^{n+1}\, dr ,
\end{multline}
respectively. Any of these functions can thus be expressed in the
form
\begin{equation}\label{repres18}
\frac 1{|x-y|^{n-1}}g\Big(x,\frac{y-x}{|y-x|}\Big) + h(x,y),
\end{equation}
where $g : \overline \Omega \times \mathbb S^{n-1} \to \mathbb R$ is
a smooth function, and $h$ is smooth for $x \neq y$, and has bounded
derivatives in $\Omega \times \Omega \setminus \{x=y\}$. As a
consequence, by \cite[Theorem 1.29]{Mik}, if $v: \Omega \to \mathbb
R$ is Lipschitz continuous, then the function $w : \Omega \to \R$
given by
 \begin{equation}\label{1.8}
 w(x) =   \int _\Omega \Big[\frac{v (y) }{|x-y|^{n-1}}g\Big(x,\frac{y-x}{|y-x|}\Big) + v(y) h(x,y)\Big]\, dy \quad \hbox{for $x \in \Omega$},
\end{equation}
 belongs to  $W^{1,1} (\Omega )$, and, for $h=1, \dots , n$, there exists a
 constant $C=C(R,h,n)$ such that
\begin{equation}\label{1.9}
\frac{\partial w}{\partial x_h}(x) =  \int _{\Omega} \frac
1{|x-y|^{n}}f\Big(x,\frac{y-x}{|y-x|}\Big) v(y) dy + C v(x) \quad
\hbox{for a.e. $x \in \Omega$,}
\end{equation}
where $f:\Omega\times \mathbb S^{n-1}\rightarrow\mathbb R$ obeys
\begin{equation}\label{repres19}
\frac 1{|x-y|^{n}}f\Big(x,\frac{y-x}{|y-x|}\Big) =
\bigg[\frac{\partial}{\partial x_h} \bigg(\frac
1{|x-y|^{n-1}}g\Big(z,\frac{y-x}{|y-x|}\Big)\bigg)\bigg]_{\lfloor
z=x}\quad\text{for }(x,y)\in\Omega\times\Omega,\,\,x\neq y.
\end{equation}
Define $N:\Omega\times(\rn\setminus \{0\})\rightarrow\mathbb R$ as
\begin{equation}\label{N}
N(x, z)= \frac 1{|z|^{n}}f\Big(x,\frac{z}{|z|}\Big)\quad\text{for
}(x,z)\in\Omega\times(\rn\setminus \{0\}).
\end{equation}
We claim that such a function fulfills properties
\eqref{N1}--\eqref{N6}. Properties \eqref{N1} and \eqref{N4} hold
trivially. Condition \eqref{N2}   holds by the results of
\cite[Section 8]{Mik}. Property \eqref{N3} is a consequence of the
smoothness of $f$. Conditions \eqref{N5} and \eqref{N6} can be shown
via standard arguments.
\par\noindent
Altogether, we have shown that equations \eqref{repres} and
\eqref{july1} hold if $\bfu \in C^\infty (\Omega , \rn)$.  We claim
that these equations continue to hold even if  $\bfu \in
E^{D,A}(\Omega, \rn)$.
Since the function $A(t)$ dominates the function $t\log(1+t)$ near
infinity, $E^{D,A}(\Omega, \rn) \to E^{D}L\log L (\Omega,\rn)$, and
hence $\bfu \in E^{D} L\log L (\Omega,\rn)$. Inasmuch as the
function $t \log (1+ t)$ satisfies the $\Delta _2$ condition near
infinity, a standard convolution argument, as, for instance, in the
proof of \cite[Proposition 1.3, Chapter 1]{Te}, tells us that
$C^\infty (\Omega , \rn)$ is dense in $ E^{D} L\log L (\Omega
,\rn)$. Thus, there exists a sequence $\{\bfu _m\} \subset C^\infty
(\Omega , \rn)$ such that
$$\bfu _m \to \bfu \quad \hbox{in $E^{D} L\log L (\Omega ,\rn)$.}$$
In particular, \begin{equation}\label{july10}
 \bfu _m \to \bfu \quad
\hbox{in $L\log L (\Omega,\rn )$,} \end{equation}
 and
\begin{equation}\label{july11} \mathcal E ^D \bfu _m \to \mathcal E^D \bfu \quad \hbox{in $L\log L (\Omega,\rn)$.}
\end{equation}
We already know that formulas \eqref{repres} and \eqref{july1} hold
with $\bfu$ replaced by $\bfu _m$. By \eqref{repres18}, all kernels
$R_{kj}$ appearing in \eqref{repres2} admit a bound of the form
\begin{align}\label{2701}
|R_{kj}(x,y)| \leq \frac C{|x-y|^{n-1}} \quad \hbox{for $x \neq
y$.}
\end{align}
 Now, recall that any  integral operator, with kernel bounded by
$\tfrac {C}{|x-y|^{n-1}}$, is  bounded  in any Orlicz space
$L^A(\Omega)$, and in particular in $ L\log L (\Omega )$. Thus, by
equations
 \eqref{july10} and \eqref{july11},  passing to the limit (possibly for a subsequence) in the representation
 formula \eqref{repres} applied to $\bfu _m$, implies that it continues to
 hold also for $\bfu$.\\
 Moreover, owing to \cite[Theorem 3.8]{BrCi}, singular integral operators whose kernel $N$
 satisfies \eqref{N1}--\eqref{N6} are bounded from $ L\log L (\Omega
 )$ into $L^1(\Omega)$. Thus,  passing to the limit in \eqref{july1} applied
 to $\bfu _m$, and  making use of \eqref{july10} and
 \eqref{july11} again, tell us that \eqref{july1} holds for $\bfu$
 as well. \qed

The proof of Lemma \ref{febbraio100} below relies upon the following
characterization of Hardy type inequalities in Orlicz spaces from
\cite{Ci2} (see also \cite{Cianchisharp, Cianchistrong} for
alternative versions).

\begin{lemma}\label{hardy} \textrm{(\cite[Lemma 5.2]{Ci2})}
 Let $A$ be and  $B$ be  Young functions, and let $L\in (0, \infty
)$.
\par\noindent
(i) There exists a constant $C$ such that
\begin{equation}\label{hardy3norm}
\bigg\|\frac 1s\int _0^s f (r)\, dr \bigg\|_{L^B(0, L)} \leq C \|f
\|_{L^A(0, L)}
\end{equation}
for every $f \in L^A(0, L)$ if and only if either $L<\infty$ and
condition \eqref{1.1} holds for some  $t_0 \geq 0$, or $L=\infty$
and \eqref{1.1} holds with  $t_0 = 0$. In particular, in the latter
case, the constant $C$ in \eqref{hardy3norm} depends only on the
constant $c$ appearing in \eqref{1.1}.
\par\noindent
 (ii) There exists a constant $C$ such that
\begin{equation}\label{hardy5young}
\bigg\|\int _s^L f (r)\, \frac{dr}r \bigg\|_{L^B(0, L)} \leq C \|f
\|_{L^A(0, L)}
\end{equation}
for every $f \in L^A(0, L)$ if and only if either $L<\infty$ and
condition \eqref{1.2} holds for some  $t_0 \geq 0$, or $L=\infty$
and \eqref{1.2} holds with  $t_0 = 0$. In particular, in the latter
case, the constant $C$ in \eqref{hardy5young} depends only on the
constant $c$ appearing in \eqref{1.2}.
\end{lemma}

\begin{lemma}\label{febbraio100}
Let $\Omega$ be a bounded open set in $\rn$, $n \geq 3$, which is
starshaped with respect to a ball. Let $A$ and $B$ be Young
functions satisfying conditions \eqref{1.1} and \eqref{1.2}.  Then
$E^{D,A}(\Omega, \rn) \subset W^{1,B}(\Omega , \rn)$. Moreover, on
denoting by $P$   the operator defined as in \eqref{repres1}, there
exists a constant $C$ such that
\begin{equation}\label{febbraio101}
\|\bfu - P\bfu \|_{L^A(\Omega ,\R^{n})} + \|\nabla (\bfu - P\bfu
)\|_{L^B(\Omega,\R^{n\times n})} \leq C \|\mathcal E^D
\bfu\|_{L^A(\Omega,\R^{n\times n})}
\end{equation}
for every $\bfu \in E^{D,A}(\Omega, \rn)$.
\end{lemma}
\par\noindent
{\bf Proof}. Let us denote by $T$ the operator defined by
\begin{equation*}
T\psi(x) = \sum _{k,j} \int _{\Omega} K_{ijhk}(x,y) \psi(y) dy \quad
\hbox{for a.e. $x \in \Omega$,}
\end{equation*}
for $\psi \in L \log L (\Omega)$, where $K_{ijhk}$ is as in Lemma
\ref{repr}. One can deduce from \cite[Theorem 3.8]{BrCi}  that there
exists a constant $C$ depending on $n$, the diameter of $\Omega$ and
the constants appearing in \eqref{N3}--\eqref{N6} such that
\begin{equation}\label{rearrestim}
(T\psi)^*(s) \leq C \bigg(\frac 1s\int _0^s \psi^*(r)\, dr + \int
_s^{|\Omega|}\psi^*(r)\, \frac {dr}r \bigg) \quad \hbox{for $s \in
(0, |\Omega|)$.}
\end{equation}
Hence, owing to \eqref{july1} and Lemma \ref{hardy}, there exists a
constant $C$ such that
\begin{equation}\label{bound2}
\|\nabla (\bfu - P\bfu )\|_{L^B(\Omega,\R^{n\times n})} \leq C \|\mathcal E^D
\bfu\|_{L^A(\Omega,\R^{n\times n})}
\end{equation}
for every $\bfu \in E^{D,A}(\Omega, \rn)$, where $P\bfu$ is defined
as in \eqref{repres1}.
\\ On the other hand, \eqref{repres} and \eqref{2701}, and the fact that  any  integral operator, with kernel bounded by
$\tfrac {C}{|x-y|^{n-1}}$, is  bounded  in any Orlicz space
$L^A(\Omega)$, ensure that
\begin{equation}\label{bound1}
\|\bfu - P\bfu \|_{L^A(\Omega ,\R^{n})} \leq C \|\mathcal E^D
\bfu\|_{L^A(\Omega,\R^{n\times n})}
\end{equation}
for some constant $C$ and every $\bfu \in E^{D,A}(\Omega, \rn)$.
Inequality \eqref{febbraio101} follows from \eqref{bound2} and
\eqref{bound1}. \qed

\begin{lemma}\label{lemma3z}
  Let $\Omega$ be a bounded connected open set with the cone property in $\rn$, $n \geq 3$, and let~$A$ be
  a Young function. Let~$\Pi : L^1(\Omega , \rn) \to \Sigma$ be a linear
  projection operator
   such that
  \begin{equation}\label{lemma3z.1}
    \|\Pi \bfu\|_{L^1(\Omega,\R^{n})}
    \leq C \|\bfu \|_{L^1(\Omega,\R^{n} )}
  \end{equation}
  for some constant~$C$, and every $\bfu \in L^1(\Omega , \rn)$. Then there exists a constant~$C'$ such that
  \begin{equation}\label{april40}
    \inf_{\bfw \in \Sigma} \|\bfu - \bfw\|_{L^A(\Omega,\R^{n})} \leq
    \|\bfu - \Pi \bfu\|_{L^A(\Omega,\R^{n})} \leq C'\, \inf_{\bfw \in
      \Sigma} \|\bfu - \bfw\|_{L^A(\Omega,\R^{n})}
\end{equation}
for every $\bfu \in L^A(\Omega,\R^{n})$, and
    \begin{equation}\label{april41}
    \inf_{\bfw \in \Sigma} \|\nabla(\bfu -
    \bfw)\|_{L^A(\Omega,\rmat)} \leq \|\nabla(\bfu - \Pi
    \bfu)\|_{L^A(\Omega,\rmat)} \leq C'\, \inf_{\bfw \in \Sigma}
    \|\nabla(\bfu - \bfw)\|_{L^A(\Omega,\rmat)}
  \end{equation}
for every $\bfu \in W^{1,A}(\Omega,\R^{n})$.
\end{lemma}
\par\noindent
{\bf Proof}. The left-wing inequalities in \eqref{april40} and
\eqref{april41} are trivial. As far as the right-wing inequalities
are concerned, given any $\bfw \in \Sigma$, and any $\bfu$ in
$L^A(\Omega,\R^{n})$, or  in $W^{1,A}(\Omega,\R^{n})$, according to
whether \eqref{april40} or \eqref{april41} is in question, set
$$\bfv =
\bfw +
  (\bfu-\bfw)_\Omega.$$ 
Here, $(\bfu-\bfw)_\Omega$ denotes the mean-value of a the vector-valued function $\bfu-\bfw$ over the set $\Omega$.
Since $\Pi$, restricted to $\Sigma$, agrees with the identity map,
  have that $\Pi \bfv = \bfv$. As a consequence,
  \begin{align*}
    \bfu - \Pi \bfu = (\bfu- \bfv) - \Pi (\bfu-\bfv).
  \end{align*}
  Thus,
  \begin{equation}\label{april42}
    \|\bfu - \Pi \bfu\|_{L^A(\Omega,\R^{n})} \leq \norm{\bfu - \bfv
    }_{L^A(\Omega,\R^{n})} + \norm{\Pi(\bfu -
      \bfv)}_{L^A(\Omega,\R^{n})},
    \end{equation}
    and
    \begin{equation}\label{april43}
    \|\nabla(\bfu - \Pi \bfu)\|_{L^A(\Omega,\rmat)} \leq
    \norm{\nabla(\bfu - \bfv)}_{L^A(\Omega,\rmat)} + \norm{\nabla
      \Pi(\bfu - \bfv)}_{L^A(\Omega,\rmat)}.
  \end{equation}
 By the triangle inequality,
\begin{equation}\label{april45}
 \norm{\bfu - \bfv
    }_{L^A(\Omega,\R^{n})} = \norm{\bfu - \bfw -
  (\bfu-\bfw)_\Omega
    }_{L^A(\Omega,\R^{n})} \leq 2 \norm{\bfu - \bfw
    }_{L^A(\Omega,\R^{n})}.
    \end{equation}
    Also,
\begin{equation}\label{april46}
 \norm{\nabla(\bfu - \bfv)}_{L^A(\Omega,\rmat)} =  \norm{\nabla(\bfu -
 \bfw)}_{L^A(\Omega,\rmat)}.
 \end{equation}
  Since the range of $\Pi$ is a finite dimensional space, where  all
  norms are equivalent, there exists a constant $C''$ such that
  \begin{align}\label{aprilC}
    \norm{\Pi(\bfu - \bfv)}_{L^A(\Omega,\R^{n})} + \norm{\nabla
      \Pi(\bfu - \bfv)}_{L^A(\Omega,\rmat)} \leq C''\, \norm{\Pi(\bfu -
      \bfv)}_{L^1(\Omega,\R^{n})}.
  \end{align}
Inequality \eqref{lemma3z.1} ensures that
 \begin{align}\label{aprilD}
    \norm{\Pi(\bfu - \bfv)}_{L^1(\Omega,\R^{n})} \leq C\, \|\bfu -
    \bfv \|_{L^1(\Omega,\R^{n} )} = C\, \norm{ \bfu-\bfw - (\bfu
        - \bfw)_\Omega}_{L^1(\Omega,\R^{n})}.
  \end{align}
  Now, by the triangle inequality,
\begin{align}\label{aprilE}
\norm{ \bfu-\bfw - (\bfu -
        \bfw)_\Omega}_{L^1(\Omega,\R^{n})} \leq 2 \norm{\bfu -
      \bfw}_{L^1(\Omega,\R^{n})}.
\end{align}
On the other hand, our assumptions on $\Omega$ ensure that a
Poincar\'e type inequality holds in $W^{1,1}(\Omega, \rn)$, and
hence there exists a constant $C$ such that
 \begin{align}\label{aprilF}
    \norm{ \bfu-\bfw - (\bfu -
        \bfw )_\Omega}_{L^1(\Omega,\R^{n})} \leq C \norm{\nabla (\bfu
      - \bfw)}_{L^1(\Omega,\rmat)}.
     \end{align}
Altogether, inequalities \eqref{april40} and \eqref{april41} follow.
\qed

\bigskip
\par\noindent

Let $A$ and $B$ be Young functions. An open set $\Omega$ in $\rn$,
$n \geq3$, will be called admissible with respect to the couple
$(A,B)$ if there exists a constant $C$ such that
\begin{equation}\label{lemma4.2}
  \inf_{\bfw \in \Sigma} \|\bfu - \bfw\|_{L^A(\Omega,\R^{n})} +
  \inf_{\bfw \in \Sigma}  \|\nabla (\bfu - \bfw)\|_{L^B(\Omega
    ,\R^{n\times n})} \leq C  \| \mathcal E^{D}
  \bfu\|_{L^A(\Omega,\R^{n\times n})}
\end{equation}
for  every $\bfu \in E^{D,A}(\Omega , \rn)$.

\smallskip

\begin{lemma}\label{lemma4}
Let $A$ and $B$ be Young functions, and let $\Omega _1$ and $\Omega
_2$ be bounded connected open sets with the cone property in $\rn$,
$n \geq 3$. Assume that each of them is admissible with respect to
$(A,B)$, and $\Omega _1 \cap \Omega _2 \neq \emptyset$. Then the set
$\Omega _1 \cup \Omega _2$ is admissible with respect to $(A,B)$ as
well.
%
%
\end{lemma}
\par\noindent
{\bf Proof}. Let $\mathcal{B}\subset \Omega_1 \cap \Omega_2$ be a
ball. Fix~$\omega \in C^\infty_0(\mathcal{B})$.  Denote by
$\mathcal{P}_2$ the space of polynomials of degree not exceeding
$2$, and by $\Pi_3 \bfu \in \mathcal{P}_2$ the averaged Taylor
polynomial of third-order with respect to~$\omega$ of a function
$\bfu \in L^1(\Omega_1 \cup \Omega_2, \rn)$ -- see~\cite{BreSco94}.
The operator $\Pi_3 : L^1(\Omega_1 \cup \Omega_2, \rn) \to
\mathcal{P}_2$ is linear, and, by \cite[Corollary~4.1.5]{BreSco94},
there exists a constant $C$ such that
\begin{align*}
  \norm{\Pi_3 \bfu}_{L^1(\Omega_1 \cup \Omega_2, \rn)} &\le C\,
  \norm{\bfu}_{L^1(\mathcal{B}, \rn)}
\end{align*}
for every $\bfu  \in L^1(\Omega_1 \cup \Omega_2, \rn)$. Furthermore,
on denoting by $\Pi_\Sigma$  the $L^2$-orthogonal projection from
$\mathcal{P}_2$ into~$\Sigma$, one has that
\begin{align*}
  \norm{\Pi_\Sigma \bfp}_{L^1(\Omega_1 \cup \Omega_2, \rn)} &\le c\,
  \norm{\bfp}_{L^1(\Omega_1 \cup \Omega_2, \rn)}
\end{align*}
for every $\bfp \in  \mathcal{P}_2$. Thus, the linear operator $\Pi
= \Pi_\Sigma \circ \Pi_3$  maps  $L^1(\Omega_1 \cup \Omega_2)$
into~$\Sigma$, and there exists a constant $C$ such that
\begin{align}\label{april47}
\norm{\Pi \bfu}_{L^1(\Omega_j, \rn)} \leq  \norm{\Pi
\bfu}_{L^1(\Omega_1 \cup \Omega_2, \rn)} \le C\,
  \norm{\bfu}_{L^1(\mathcal{B})} \leq C\,
  \norm{\bfu}_{L^1(\Omega_j, \rn)} \quad j=1,2
\end{align}
for every $\bfu \in L^1(\Omega_1 \cup \Omega_2, \rn)$. Owing to
inequality \eqref{april47},  Lemma~\ref{lemma3z} ensures that there
exists a constant $C$ such that
\begin{align}\label{april48}
  \inf_{\bfw \in \Sigma} \norm{\bfu - \bfw}_{L^A(\Omega_1 \cup
    \Omega_2, \rn)} &\leq \norm{\bfu - \Pi \bfu}_{L^A(\Omega_1 \cup
    \Omega_2, \rn)} \leq \sum_{j=1,2} \norm{\bfu - \Pi \bfu}_{L^A(\Omega_j, \rn)}
\\ \nonumber &
    \leq C\,\sum_{j=1,2}\, \inf_{\bfw \in \Sigma} \norm{\bfu
    - \bfw}_{L^A(\Omega_j, \rn)}
    \end{align}
     for every $\bfu \in L^A(\Omega, \rn)$. Similarly, by
     Lemma~\ref{lemma3z} applied with $A$ replaced by $B$, there
     exists a constant $C$ such that
 \begin{align}\label{april49}
  \inf_{\bfw \in \Sigma} \norm{\nabla(\bfu - \bfw)}_{L^B(\Omega_1 \cup
    \Omega_2, \rmat)} &\leq  \norm{\nabla(\bfu - \Pi
    \bfu)}_{L^B(\Omega_1 \cup \Omega_2, \rmat)}
 \leq \sum_{j=1,2} \norm{\nabla(\bfu - \Pi
    \bfu)}_{L^B(\Omega_j, \rmat)}
    \\ \nonumber &
    \leq C\,\sum_{j=1,2} \,\inf_{\bfw
    \in \Sigma} \norm{\nabla(\bfu - \bfw)}_{L^B(\Omega_j, \rmat)}.
\end{align}
The conclusion follows from \eqref{april48}--\eqref{april49}, and
\eqref{lemma4.2} applied with $\Omega = \Omega _j$, for $j=1,2$.
\qed

\begin{lemma}\label{suff}
  Let $\Omega$ be a connected, bounded open set with the cone
  property in $\rn$, $n \geq 3$.  Let $A$ and $B$ be Young functions satisfying \eqref{1.1}
  and \eqref{1.2}. Then~$\Omega$ is admissible with respect
  to~$(A,B)$. Moreover, if $\Pi : E^{D,A}(\Omega , \rn) \to \Sigma$
  is a linear projection operator such that
  \begin{equation}\label{lemma3.1}
    \|\Pi \bfu\|_{L^1(\Omega,\R^{n})}
    \leq C \|\bfu \|_{L^1(\Omega,\R^{n} )}
  \end{equation}
  for some constant $C$ and every $\bfu \in E^{D,A}(\Omega
  ,\R^{n})$, then there exists a constant $C'$ such that
  \begin{equation}\label{suff2}
    \|\bfu - \Pi \bfu\|_{L^A(\Omega
      ,\R^{n})} +     \|\nabla (\bfu - \Pi \bfu)\|_{L^B(\Omega,\R^{n\times n}
      )} \leq C' \|\mathcal E^D \bfu\|_{L^A(\Omega,\R^{n\times n})}
  \end{equation}
  for every $\bfu \in E^{D,A}(\Omega, \rn)$.  In particular,
  inequality \eqref{suff2} holds with $\Pi =P$, where $P$ is  defined as in \eqref{repres}--\eqref{repres1}.
\end{lemma}
\par\noindent
{\bf Proof}. The statement holds if $\Omega $ is starshaped with
respect to a ball, thanks Lemma~\ref{febbraio100} and
Lemma~\ref{lemma3z}, applied with $\Pi =P$.
On the other hand, any open set $\Omega$ as in the statement is the
finite union of open sets $\Omega_i$, $i=1, \dots , k$, starshaped
with respect to a ball. Since $\Omega$ is connected, after,
possibly, relabeling, we may assume that, the sets $\cup
_{i=1}^{j-1}\Omega _i$ and $\Omega _j$ have a non-empty
intersection. The conclusion then follows from repeated use of
Lemma~\ref{lemma4}. \qed

\medskip
\noindent {\bf Proof of Theorem \ref{sufficiency0}}.
Let $\mathcal{B}'$ be an open  ball such that~$\overline \Omega
\subset \mathcal{B}'$. Let $2\mathcal{B}'$ denote the ball with same
center as $\mathcal{B}'$, and twice its radius.  Since $\bfu \in
E^{D,A}_0(\Omega, \rn)$, its extension by zero to~$2\mathcal{B}'$,
still denoted by $\bfu$, belongs to $E^{D,A}_0(2\mathcal{B}', \rn)$.
 Let $\mathcal{B}$ be a ball in $2\mathcal{B}' \setminus
\mathcal{B}'$, and pick any function $\omega \in
C^\infty_0(\mathcal{B})$. Let $\Pi=\Pi_\Sigma \circ \Pi_3$ be the
projection operator defined as in the proof of Lemma~\ref{lemma4}.
In particular, $\Pi : L^1(2\mathcal{B}', \rn) \to \Sigma$, and
\begin{align*}
  \norm{\Pi \bfu}_{L^1(2\mathcal{B}', \rn)} &\le C\,
  \norm{\bfu}_{L^1(\mathcal{B}, \rn)}
\end{align*}
for some constant $C$. Hence, since $\bfu =0$ in~$\mathcal{B}$, we
infer that $\Pi \bfu=0$. The conclusion is now a consequence of
Lemma~\ref{suff}. \qed



\medskip
\par\noindent
{\bf Proof of Theorem \ref{sufficiency}}. 
The conclusion follows from Lemmas ~\ref{suff} and \ref{lemma3z}.
\qed

\medskip
\par\noindent

As a byproduct of our approach to Theorems \ref{sufficiency0} and
\ref{sufficiency}, one can derive the  Poincar\'e type inequalities
in $E^{D,A}_0(\Omega,\rn)$ and $E^{D,A}(\Omega,\rn)$, of independent
interest, which are stated in the next theorem. Let us emphasize
that
 they hold for any Young function $A$. The special case when $A(t) =
 t$ was considered in \cite{FR}.

\begin{theorem}\label{poincareineq}
  Let $\Omega $ be a bounded open set in $\rn$, $n \geq 3$, and let $A$ be any Young
  function. Then there exists a constant $C$ such that
      \begin{equation}
      \label{poincareineq2} \|\bfu\|_{L^A(\Omega
        ,\R^{n})}  \leq C \|\mathcal E^D \bfu\|_{L^A(\Omega,\R^{n\times n})}
    \end{equation}
 for every $\bfu \in E^{D,A}_0(\Omega,\rn)$.
\\ Assume in addition that $\Omega$ is connected and has the cone
property. Then there exists a constant $C$ such that
\begin{equation}\label{poincareineq1}
  \inf_{\bfw \in \Sigma} \|\bfu - \bfw\|_{L^A(\Omega,\R^{n})} \leq C  \| \mathcal E^{D}
  \bfu\|_{L^A(\Omega,\R^{n\times n})}
\end{equation}
for  every $\bfu \in E^{D,A}(\Omega , \rn)$.
\end{theorem}

\medskip
\par\noindent
{\bf Proof, sketched}. A proof of inequalities \eqref{poincareineq2}
and \eqref{poincareineq1} can be accomplished along the same lines
as the proof of inequalities \eqref{1.3.0} and \eqref{1.3main}. One
has just to make use of the inequalities in the statement of Lemmas
\ref{febbraio100}--\ref{suff}, and in the definition of admissible
domains, without  gradient norms on their left-hand sides. Since
inequality \eqref{bound1} does not require any assumption on $A$,
the relevant lemmas, and hence Theorem \ref{poincareineq}, hold for
any Young function $A$. The details are omitted for brevity. \qed

\medskip
\par\noindent

\section{Necessary  conditions  for  Korn-type and  related  inequalities}
\label{sec:Korn2}

 The key step in the proof of the necessity of assumptions \eqref{1.1} and
 \eqref{1.2} in Theorems \ref{thm:main0} and \ref{thm:main}, as well
 as in the other statements  of Section \ref{main}, is the following
 results, dealing with Korn-type inequalities for functions
 subject to vanishing boundary conditions.

\begin{theorem}\label{necessity}
Let $\Omega $ be a  bounded open set in $\rn$, $n \geq 2$. Let $A$
and $B$ be Young functions such that
\begin{equation}\label{1.3nec}
\| \nabla \bfu \|_{L^B(\Omega , \rmat)} \leq C \| \mathcal E
\bfu\|_{L^A(\Omega , \rmat)}
\end{equation}
for some constant $C$, and for every $\bfu \in W^{1,1}_0(\Omega ,
\rn) \cap E^{1,A}_0(\Omega , \rn)$.  Then conditions \eqref{1.1} and
\eqref{1.2} hold.
\end{theorem}

\smallskip
\par
A proof of inequalities \eqref{1.1} and \eqref{1.2}  rests upon
different choices of trial functions in inequality \eqref{1.3nec}.
In particular, our derivation of \eqref{1.1} is related to an
argument from \cite{CFM}, which makes use of the so called \lq\lq
laminates" to provide an alternative proof of the failure of the
Korn inequality \eqref{intro1} for $p=1$.\\ A first-order laminate
is a probability measure $\nu$ on $\R^{n\times n}$ of the form
$$\nu=\lambda\delta_\bfA+(1-\lambda)\delta_\bfB,$$
where $\lambda\in(0,1)$, $\bfA , \bfB \in \rmat$, and
$\rank(\bfA-\bfB)=1$. Here $\delta_\bfX$ denotes the Dirac measure
on $\R^{n\times n}$ concentrated at the matrix $\bfX$. The matrix
$\lambda \bfA+(1-\lambda)\bfB$ is called the average of the laminate
$\nu$. A second-order laminate is obtained on replacing
$\delta_\bfA$ [resp. $\delta_\bfB$] with a first-order laminate with
average $\bfA$ [$\bfB$]. Higher-order laminates are defined
accordingly via an iteration process. We refer to \cite{KMS} and
\cite{M} for a detailed discussion on laminates. The following
approximation lemma from \cite{CFM} will be exploited in our proof
of Theorem \ref{necessity}.

\begin{lemma}
  \label{lems2}  {\rm (\cite[Equation (5)]{CFM})} Let $\nu$ be a laminate in $\R^{n\times n}$ with average $C$, and let $r>0$ Then there
  exists a sequence $\{\bfu _i\}$ of uniformly Lipschitz continuous functions
  $\bfu_i: (0,r)^n\rightarrow \R^n$, such that $\bfu_i (x)= Cx$ for $x \in \partial (0,r)^n$, and
  \begin{equation}\label{feb2}
  \lim _{i \to \infty} \int_{(0,r)^n}\Phi(|\nabla \bfu_i|)\, dx = r^{n}\int_{\R^{n\times
  n}}\Phi(|\bfX|)\, d\nu(\bfX),\end{equation}
  for every continuous function $\Phi$.
\end{lemma}

\medskip
\par\noindent {\bf Proof of Theorem \ref{necessity}}. \emph{Part 1}: Inequality \eqref{1.2} holds.
%
%
\par\noindent
Assume, without loss of generality, that the unit ball ${\mathcal
B}_1$, centered at $0$, is contained in $\Omega$, and denote by
$\omega _n$ its Lebesgue measure.  Let us preliminarily observe that
inequality \eqref{1.3nec} implies that
\begin{equation}\label{A>B}
\hbox{$A$ dominates $B$ near infinity.}
\end{equation}
 Indeed, given any  nonnegative function
$h\in L^A (0, \omega _n)$, consider the function $\bfv : {\mathcal
B}_1 \to \rn$ given by
$$\bfv (x) = \bigg(\int _{|x|}^1 h(\omega _n r^n)\,
dr\,, 0, \dots , 0\bigg) \quad \hbox{for $x \in {\mathcal B}_1$.}$$
Then $\bfv \in L^A({\mathcal B}_1, \rn)$, and
$$|\mathcal E\bfv (x)| \leq |\nabla \bfv (x)| =h(\omega _n |x|^n) \qquad \hbox{for $x \in {\mathcal B}_1$.}$$
An application of \eqref{1.3nec}, with $\bfu$ replaced by $\bfv$,
thus tells us that
$$\|h\|_{L^B (0, \omega _n)} = \| \nabla \bfv \|_{L^B(\Omega , \rmat)} \leq C \| \mathcal E
\bfv\|_{L^A(\Omega , \rmat)} \leq C \| \nabla \bfv\|_{L^A(\Omega ,
\rmat)} = \|h\|_{L^A (0, \omega _n)}.$$ Thus $L^A(0, \omega _n) \to
L^B(0, \omega _n)$, and  \eqref{A>B} follows.
\\
Now, given  $h$ as above, define the function $\rho : [0, 1] \to [0,
\infty]$ as
\begin{align*}
\rho(r)=\int_r^1\frac{h(\omega _nt^n)}{t}\dt \quad \hbox{for $r \in
[0, 1]$,}
\end{align*}
and the function  $\bfu : {\mathcal B}_1 \to \rn$ as
\begin{align*}
  \bfu(x) = \bfQ \,x \rho(\abs{x}) \quad \hbox{for $x \in {\mathcal B}_1$,}
\end{align*}
where  $\bfQ \in \setR^{n\times n}$ is any skew-symmetric matrix
such that $\abs{\bfQ}=1$. One has that $\bfu$ is a weakly
differentiable function, and
\begin{align*}
  \mathcal E\bfu(x) &=\frac{\bfQ x\otimes ^{\rm sym} x}{|x|^2}\rho'(\abs{x})\abs{x},
  \\
  \nabla \bfu(x) &= \bfQ \rho(\abs{x})+\frac{\bfQ x\otimes x}{|x|^2}\rho'(\abs{x})\abs{x}
\end{align*}
for a.e. $x \in {\mathcal B}_1$. Here, $\otimes ^{\rm sym}$ denotes
the symmetric part of the tensor product of two vectors in $\rn$.
Hence,
\begin{align*}
  \abs{\mathcal E \bfu(x)} &\leq|\rho'(\abs{x})|\abs{x} = h(\omega _n |x|),
  \\
  \rho(\abs{x}) & \leq  \abs{\nabla \bfu (x)} +
  |\rho'(\abs{x})|\abs{x}= \abs{\nabla \bfu (x)} + h(\omega _n |x|)
\end{align*}
for a.e. $x \in {\mathcal B}_1$. Thus, owing to  \eqref{1.3nec} and
\eqref{A>B},
\begin{align}\label{feb1}
&\Big\|\int_{s}^{\omega _n}\frac{h(r)}{r}\,dr\Big\|_{L^B(0, \omega
_n)}= \Big\|\int_{|x|}^{1}\frac{h(\omega _n
t^n)}{t}\,dt\Big\|_{L^B({\mathcal B}_1)} = \|\rho (|x|)\|_{L^B({\mathcal B}_1)} \\
\nonumber &\quad \leq \|\nabla \bfu \|_{L^B({\mathcal B}_1,\rmat)} + \|h(\omega
_n |x|^n) \|_{L^B({\mathcal B}_1)} \leq C \|\mathcal E \bfu
\|_{L^A({\mathcal B}_1,\rmat)} + \|h(\omega _n |x|^n) \|_{L^A({\mathcal B}_1)}\\
\nonumber &\quad \leq C'\|h(\omega _n |x|^n) \|_{L^A({\mathcal B}_1)} =
C' \|h(s) \|_{L^A(0, \omega _n)}
\end{align}
for suitable constants $C$ and $C'$. Thanks to the arbitrariness of
$h$, inequality \eqref{feb1} implies, via Lemma \ref{hardy}, that
\eqref{1.2} holds for some $c$ and $t_0$.

\smallskip
\par\noindent
\emph{Part 2}: Inequality \eqref{1.1} holds.
\\ Let us preliminarily  note that, if $A(t)=\infty$ for
large $t$, then \eqref{1.1} holds trivially. We may thus assume that
$A$ is finite-valued, and  hence  continuous. By \eqref{A>B},  the
function $B$ is also finite-valued and  continuous.
\par\noindent For
ease of notations, we hereafter focus on
   case when $n=2$.  An analogous argument carries over to any dimension along the lines of
  \cite[Lemma 3]{CFM}.  Given $a, b \in \mathbb R$,
  define the matrix $\bfG_{a,b}$ as
  $$\bfG_{a,b} = \left(\begin{matrix} 0 & a \\
  b & 0 \end{matrix}\right),$$
  and set
  $\delta_{a,b} =
  \delta_{\bfG_{a,b}}$. Next, define the sequence $\{\mu^{(m)}\}$ of
  laminates of order $2m$ by iteration as:
  \begin{equation} \begin{cases}
    \mu^{(0)} & = \delta_{t,t},
    \\
    \mu^{(m)} & =
    \tfrac{1}{3}\delta_{2^{-m}t,-2^{-m}t}+ \tfrac{1}{6}
    \delta_{-2^{1-m}t,2^{1-m}t}+\tfrac{1}{2} \mu^{(m-1)}
    \end{cases}
  \end{equation}
  for $m \in \setN$. We claim that
 $\mu^{(m)}$ is a
  laminate with average $\bfG_{2^{-m}t,2^{-m}t}$ for $m \in \mathbb N$.
Indeed, one has that \begin{equation}\label{telc1}\mu^{(m)}  =
    \tfrac{1}{4} \delta_{-2^{1-m}t,2^{1-m}t}+\tfrac{3}{4}
    \mu^{(m-1)}.\end{equation}
     Since  $\rank(\bfG_{-t,t}-\bfG_{t,t})=1$,
  the right-hand side of \eqref{telc1} is  a laminate with average
  $\bfG_{2^{-1}t,t}$ for $m=1$. Hence, $\mu^{(1)}$ is a laminate with
  average
  $\bfG_{2^{-1}t,2^{-1}t}$. An induction argument then proves our claim.
  Now, note the representation formula
  \begin{align}
    \label{eq:mun}
    \mu^{(m)} &= 2^{-m} \delta_{t,t} + \sum_{k=1}^m \big( \tfrac{1}{3}
    2^{k-m} \delta_{2^{-k}t,-2^{-k}t}+ \tfrac{1}{6} 2^{k-m}
    \delta_{-2^{1-k}t,2^{1-k}t} \big)
  \end{align}
  for $m \in \setN$. Observe that $\delta_{t,t}$ is concentrated at
  a symmetric matrix, whereas the sum in \eqref{eq:mun} is
  concentrated at skew-symmetric matrices. Define the functions
  $\Phi _j : \mathbb R^{2\times 2} \to [0, \infty)$, for $j=1,2$, as
  \begin{align*}
    \Phi_1(\bfX) & = A(\abs{\bfX^\sym-\bfG_{2^{-m}t,2^{-m}t}}),
    \\
    \Phi_2(\bfX) & =
    B(C^{-1}\abs{\bfX-\bfG_{2^{-m}t,2^{-m}t}}),
  \end{align*}
for $\bfX \in \mathbb R^{2\times 2}$. Here, $\bfX^\sym = \tfrac
12(\bfX + \bfX ^T)$, the symmetric part of $\bfX$, and $C$ is the
constant appearing in \eqref{1.3nec}. Fix $m \in \mathbb N$. Without
loss of generality, we may assume that $0\in\Omega$. Choose $r>0$ so
small  that $(0,r)^2\subset\Omega$. Given any $m \in \mathbb N$,
owing to
 Lemma~\ref{lems2} applied with  $\nu = \mu^{(m)}$,  there exists
   a sequence $\{\bfu _i\}$ of Lipschitz continuous functions $\bfu_i\,:\, (0,r)^2 \to \setR^2$,  such
  that $\bfu_i (x) = \bfG_{2^{-m}t,2^{-m}t}x$ on $\partial (0,r)^2$,
  and
  \begin{align}
    \label{eq:limPhii}
    \lim_{i \to \infty} \int_{(0,r)^2}\Phi_j(\nabla \bfu_i)\dx &=r^{2}
    \int_{\R^{2\times 2}}\Phi_j(\bfX)\,d\mu^{(m)}(\bfX)  \qquad \hbox{for
    $j=1,2$.}
  \end{align}
Define the sequence $\{\bfv _i\}$ of functions $\bfv_i : \Omega \to
\mathbb R$ as  $\bfv_i(x)  = \bfu_i(x) - \bfG_{2^{-m}t,2^{-m}t}x$ if
$x \in  (0,r)^2$, and $\bfv_i(x) =0$ if $\Omega \setminus (0,r)^2$.
Then $\bfv_i
  \in W^{1,\infty}_0(\Omega)$, and, by~\eqref{eq:limPhii},
  \begin{align}
    \label{eq:lim1}
  \lim_{i \to \infty} \int_{\Omega}A(\abs{\mathcal E\bfv_i})\dx & = \lim_{i \to \infty} \int_{(0,r)^2}A(\abs{\mathcal E\bfv_i})\dx
    \\ \nonumber & =r^{2} \int_{\R^{2\times 2}}A(\abs{
      (\bfX^\sym-\bfG_{2^{-m}t,2^{-m}t})})\,d\mu^{(m)}(\bfX),
    \\
    \label{eq:lim2}
     \lim_{i \to \infty} \int_{\Omega}B(C^{-1}\abs{\nabla \bfv_i})\dx
   & =
    \lim_{i \to \infty} \int_{(0,r)^2}B(C^{-1}\abs{\nabla \bfv_i})\dx
    \\ \nonumber & =r^{2} \int_{\R^{2\times
        2}}B(C^{-1}\abs{\bfX-\bfG_{2^{-n}t,2^{-n}t}})\,d\mu^{(m)}(\bfX).
  \end{align}
  The following chain holds:
  \begin{align}\label{march10}
    \lefteqn{\int _{\R^{2\times 2}}A\big(\abs{\bfX^{\sym}-\bfG_{2^{-m}t,2^{-m}t}}\big)
      \,d\mu^{(m)}(\bfX)} \qquad &
    \\ \nonumber
    &\leq \frac{1}{2} \int _{\R^{2\times 2}}A\big(2 \abs{\bfX^{\sym}}) \,d\mu^{(m)}(\bfX) +
    \frac{1}{2} \int _{\R^{2\times 2}} A\big(2 \abs{\bfG_{2^{-m}t,2^{-m}t}})
    \,d\mu^{(m)}(\bfX)
    \\ \nonumber
    &= \frac{1}{2} 2^{-m} A(2 \abs{ \bfG_{t,t}}) + \frac{1}{2}
    A(2 \abs{ \bfG_{2^{-m}t,2^{-m}t}})
\\ \nonumber &= \frac{1}{2} 2^{-m} A(2 \abs{ \bfG_{t,t}}) + \frac{1}{2}
    A(2 \, 2^{-m}\abs{ \bfG_{t,t}})
    \\ \nonumber
    &\leq 2^{-m} A(2 \abs{ \bfG_{t,t}}),
  \end{align}
  where the first inequality holds since $A$ is convex, the first equality holds owing to \eqref{eq:mun} and to the fact
  that $\mu^{(m)}$ is a probability measure, and the last inequality
  follows from \eqref{lambdaA}.
Coupling \eqref{eq:lim1} with \eqref{march10} yields
  \begin{align}\label{march11}
    \lim_{i \to \infty} \int_{\Omega}A(\abs{\mathcal E\bfv_i})\dx
    &\leq r^{2} 2^{-m} A(2 \abs{ \bfG_{t,t}}).
  \end{align}
  Since $A$ is a continuous function, there exists  $t_m\in(0,\infty)$ such that
  \begin{align}
    \label{eq:choicer}
    r^{2} 2^{-m} A(2 \abs{  \bfG_{t_m,t_m}})= \tfrac{1}{2}.
  \end{align}
Thanks \eqref{lambdaA}, there exists $t_0>0$, independent of $m$,
such that
  \begin{align}
    \label{eq:tm}
    t_m\leq \,t_02^m.
  \end{align}
Therefore, by neglecting, if necessary, a finite number of terms of
  the sequence $\{\bfv_i\}$, we can assume that
  \begin{align*}
    \int_{\Omega}A(\abs{\mathcal E\bfv_i})\dx \leq 1
  \end{align*}
  for $i \in \setN$. Hence,   $\norm{\mathcal E\bfv_i}_{A}
  \leq 1$ for $i \in \setN$, and,  by \eqref{1.3nec},  $\norm{\nabla \bfv_i}_{B} \leq
  C$ for $i \in \setN$. Thus,
  \begin{align*}
    \int_{\Omega}B( C^{-1} \abs{\nabla \bfv_i} )\dx \leq 1
  \end{align*}
  for   $i \in \setN$. Combining the latter inequality with equation \eqref{eq:lim2}
  tells us that
  \begin{align}\label{telc2}
    r^{2} \int_{\R^{2\times 2}}B(C^{-1}\abs{\bfX-\bfG_{2^{-m}t_m,2^{-m}t_m}})\,d\mu ^{(m)}(\bfX) \leq 1.
  \end{align}
  Next, one can make use of ~\eqref{eq:mun} and derive the following
  chain:
  \begin{align}\label{eq:1902}
    r^{-2} &\geq \int _{\R^{2\times 2}}
    B(C^{-1}\abs{\bfX-\bfG_{2^{-m}t_m,2^{-m}t_m}})\,d\mu^{(m)}(\bfX)
    \\  \nonumber
    & \geq  2^{-m}B\big(C^{-1}(1-2^{-m}) \abs{\bfG_{t_m,t_m}}\big)
    +\sum_{k=1}^m  \tfrac{1}{3} 2^{k-m}
    B\big(C^{-1}(2^{-k}-2^{-m}) \abs{\bfG_{t_m,t_m}}\big)
    \\  \nonumber
    &\quad+ \sum_{k=1}^m \tfrac{1}{6} 2^{k-m} B \big(C^{-1}
    (2^{1-k}-2^{-m}) \abs{\bfG_{t_m,t_m}}\big)
\\
        \nonumber &\geq \sum_{k=1}^{m-1} \tfrac{1}{3} 2^{k-m}
    B\big(C^{-1}(2^{-k}-2^{-m}) \abs{\bfG_{t_m,t_m}}\big)
    \\  \nonumber
    &\geq \sum_{k=1}^{m-1} \tfrac{1}{3} 2^{k-m} B\big(C^{-1}
    2^{-k-1} \abs{\bfG_{t_m,t_m}}\big)\\  \nonumber
& \geq \sum_{k=1}^{m-1} \tfrac{1}{3} \frac 1{2C} 2^{-m}t_m
\frac{B\big(\tfrac{1}{2C}
    2^{-k} t_m\big)}{\tfrac{1}{2C}2^{-k}t_m}.
  \end{align}
From   \eqref{eq:choicer}, \eqref{telc2} and \eqref{eq:1902} one
infers that
\begin{align*}
2\cdot2^{-m} A(2 \abs{  \bfG_{t_m,t_m}})\geq\sum_{k=1}^{m-1}
\tfrac{1}{3}\tfrac{1}{2C}2^{-m}t_m \frac{B\big(\tfrac{1}{2C}
    2^{-k} t_m\big)}{\tfrac{1}{2C}2^{-k}t_m}.
  \end{align*}
Hence, by \eqref{eq:tm},
\begin{align}\label{march12}
\,A(c'' t_m)\geq\,c\,t_m\,\sum_{k=1}^{m-1} \frac{B\big(\tfrac{1}{2C}
    2^{-k} t_m\big)}{\tfrac{1}{2C}2^{-k}t_m}
\geq
c'\,t_m\int_{2^{-m}\frac{t_m}{4C}}^{\frac{t_m}{4C}}\frac{B(s)}{s^2}\,ds\geq
c'\,t_m\int_{\frac{t_0}{2C}}^{\frac{t_m}{4C}}\frac{B(s)}{s^2}\,ds.
  \end{align}
  for suitable positive constants $c$, $c'$, $c''$.
Since  $\lim _{m \to \infty} t_m =\infty$, one can find $\widehat t
\geq \tfrac {t_0}{2C}$ such that, if $t> \widehat t$, then there
exists $m\in\N$ such that $t_m\leq t< t_{m+1}$.  Moreover, $\widehat
t$ can be chosen so large that  $A$ is invertible on $[\widehat
t,\infty)$ and
\begin{align*}
t_m=c_1 A^{-1}(c_2 2^m)
\end{align*}
for some positive constants $c_1,c_2$. By \eqref{lambdaA-1}, the
latter equation ensures that $t_{m+1}\leq 2 t_m$ for $m \in \mathbb
N$. Thus, owing to inequality \eqref{march12},
\begin{align*}
A(2c''t)\geq A(2c'' t_m)\geq \,A(c''t_{m+1}) \geq
c'\,t_{m+1}\int_{\frac{t_0}{2C}}^{\frac{t_{m+1}}{4C}}\frac{B(s)}{s^2}\,ds\geq
c'\,t\int_{\frac{t_0}{2C}}^{\frac{t}{4C}}\frac{B(s)}{s^2}\,ds \quad
\hbox{for $t\geq \widehat t$.}
  \end{align*}
Hence, inequality
  \eqref{1.1} follows for suitable constants $c$ and $t_0$.
\qed

\medskip
\par
The next statement is a corollary of  Theorem \ref{necessity}.

\begin{corollary}\label{necessitycor}
 Let $A$
and $B$ be Young functions. Assume that any of the following
properties holds:
\par\noindent (i)
There exists a constant $C$ such that
\begin{equation}\label{july100}
\inf _{\bfv \in \mathcal R}  \| \nabla \bfu - \nabla
\bfv\|_{L^B(\Omega , \rmat)} \leq C \| \mathcal E \bfu\|_{L^A(\Omega
, \rmat)}
\end{equation}
for some bounded connected open set $\Omega $  in $\rn$, $n \geq 2$,
and every $\bfu \in W^{1,1}(\Omega, \rn) \cap E^{A}(\Omega , \rn)$.
\par\noindent (ii)
There exists a constant $C$ such that
\begin{equation}\label{july101}
  \| \nabla \bfu\|_{L^B(\Omega , \rmat)} \leq C \| \mathcal E^D \bfu\|_{L^A(\Omega
, \rmat)}
\end{equation}
 for  some bounded open set
$\Omega$ in $\rn$, $n \geq 3$, and every $\bfu \in W^{1,1}_0(\Omega
, \rn) \cap E^{D,A}_0(\Omega , \rn)$.
\par\noindent (iii)
There exists a constant $C$ such that
\begin{equation}\label{july102}
\inf _{\bfw \in \Sigma}  \| \nabla \bfu - \nabla \bfw\|_{L^B(\Omega
, \rmat)} \leq C \| \mathcal E^D \bfu\|_{L^A(\Omega , \rmat)}
\end{equation}
 for some connected open set $\Omega $  in $\rn$, $n \geq 3$, and  every $\bfu \in W^{1,1} (\Omega ,
\rn) \cap E^{D,A}(\Omega , \rn)$.
\\
Then  conditions \eqref{1.1} and \eqref{1.2} hold.
\end{corollary}
\par\noindent
{\bf Proof}. 
Assume that (ii) holds. Then the claim follows by
Theorem~\ref{necessity}, since~$\abs{\mathcal{E}^D \bfu} \leq 2
\abs{\mathcal{E} \bfu}$.
\\ Next, suppose that (iii) holds. Let $\mathcal B'$ be a ball such that $\mathcal B' \subset \subset
\Omega$. Pick a ball $\mathcal{B}$ contained in~$\Omega \setminus
\overline{\mathcal{B}'}$, fix any function $\omega \in
C^\infty_0(\mathcal{B}')$. Given any function $\bfu \in
W^{1,1}_0(\mathcal{B} , \rn) \cap E^{D,A}_0(\mathcal{B} , \rn)$, its
continuation by zero outside $\mathcal{B}$, still denoted by $\bfu$,
belongs to $W^{1,1}_0(\Omega , \rn) \cap E^{D,A}_0(\Omega , \rn)$.
Now let $\Pi=\Pi_\Sigma \circ \Pi_3$ be the projection operator
associated with $\omega$ as in the proof of Lemma~\ref{lemma4}. In
particular, $\Pi$ maps $L^1(\mathcal{B}')$ into~$\Sigma$, and there
exists a constant $C$ such that
\begin{align*}
  \norm{\Pi \bfu}_{L^1(\Omega, \rn)} &\le C\,
  \norm{\bfu}_{L^1(\mathcal{B}', \rn)}.
\end{align*}
Since $\bfu =0$ in ~$\mathcal{B}'$, one has that~$\Pi \bfu = 0$.
Thus, by Lemma~\ref{lemma3z}, property (ii) holds with~$\Omega$
replaced by $\mathcal{B}$. Hence, the conclusion follows.
\\
Finally, assume that (i) is in force. Given any function $\bfu \in
E^{A}_0(\Omega , \rn)$,  one has that $(\nabla \bfu)_\Omega=0$.
Consequently,
\begin{align*}
 \| \nabla \bfu\|_{L^B(\Omega , \rmat)} &=  \| \nabla \bfu -(\nabla \bfu)_\Omega\|_{L^B(\Omega , \rmat)}\leq \,C \inf _{\bfS\in \mathbb R^{n\times n}}  \| \nabla \bfu - \bfS\|_{L^B(\Omega , \rmat)}\\&\leq C \inf _{\bfv \in \mathcal R}  \| \nabla \bfu - \nabla
\bfv\|_{L^B(\Omega , \rmat)}\leq C'  \|
  \mathcal E \bfu\|_{L^A(\Omega, \mathbb R^{n\times n})}\,,
  \end{align*}
for some constants $C$ and $C'$. Thus inequality \eqref{1.3nec}
holds, and the conclusion follows via Theorem \ref{necessity}. \qed

\section{Proofs of the main results}\label{proofs}

With the results of Sections \ref{sec:Korn2suff} and \ref{sec:Korn2}
at our disposal, the proofs of  Theorems \ref{thm:main0} and
\ref{thm:main} can be promptly accomplished. The  necessity of
condition \eqref{1.1} and \eqref{1.2} in Theorems
\ref{thm:plainkorn}, \ref{negativenecsuf} and \ref{corbogov} also
easily follows.

\medskip
\par\noindent
{\bf Proof of Theorem \ref{thm:main0}}. Condition (i) implies (ii)
by Theorem \ref{sufficiency0}. The reverse implication holds owing
to Corollary \ref{necessitycor}, condition (ii).
\\ In order to verify that
property (iii) implies (ii), observe that, if $\bfu$ is any function
such that $\| \mathcal E^D \bfu\|_{L^A(\Omega , \rmat)} \leq 1$,
then $\int _\Omega A(|\mathcal E^D \bfu |)\, dx \leq 1$. Hence, by
inequality \eqref{1.3.1},
$$\int_\Omega B(|\nabla\bfu|/C)\,dx\leq\,C_1+1.$$
By property \eqref{lambdaA} of Young functions, this inequality
implies that $\| \nabla \bfu\|_{L^B(\Omega , \rmat)} \leq C(C_1+1)$.
Hence, inequality \eqref{1.3.0} follows.
\\  Finally, assume that (i) is in force. Suppose first that $t_0=0$ in
\eqref{1.1} and \eqref{1.2}. We already know that inequality
\eqref{1.3.0}  holds. An inspection of the proof of Theorem
\ref{sufficiency0} and of the statement of Lemma \ref{hardy} tells
us that the constant $C$ in \eqref{1.3.0} depends only on $\Omega$
and on the constant $c$ appearing in conditions \eqref{1.1} and
\eqref{1.2}. These conditions continue to hold if the functions $A$
and $B$ are replaced with the functions $A_M$ and $B_M$ given by
$A_M(t)=A(t)/M$ and $B_M(t)=B(t)/M$ for some positive constant $M$.
Given a function $\bfu \in E^A_0(\Omega , \rn)$, set
$$M= \int _\Omega  A(C|\mathcal E^D \bfu |)\, dx\,.$$

If $M= \infty$, then inequality \eqref{1.3.1}  holds trivially. We
may thus assume that
$$\|\mathcal E^D \bfu \|_{L^{A_M}(\Omega , \rmat)} \leq 1,$$
whence, by inequality \eqref{1.3.0} applied with  $A$ and $B$
replaced by $A_M$ and $B_M$, we deduce that
\begin{equation}\label{febbraio104}
\int_\Omega B(|\nabla\bfu|)\,dx \leq  \int _\Omega  A(C| \mathcal
E^D \bfu |)\, dx\,,
\end{equation}
 namely \eqref{1.3.1}, with $C_1=0$.
 \\ Assume next that \eqref{1.1} and \eqref{1.2} just hold for some $t_0
 >0$. The functions $A$ and $B$ can be replaced with new Young functions
 $\overline A$ and $\overline B$, equivalent to $A$ and $B$ near infinity, and such
 that \eqref{1.1} and \eqref{1.2} hold for the new functions with
 $t_0=0$. The same argument as above yields \eqref{febbraio104} with
$A$ and $B$ replaced with
 $\overline A$ and $\overline B$, namely
\begin{equation}\label{febbraio105}
\int_\Omega \overline B(|\nabla\bfu|)\,dx \leq  \int _\Omega
\overline A(C| \mathcal E^D \bfu |)\, dx\,
\end{equation}
for some constant $C$. Since $\overline A$ and $\overline B$ are
equivalent to $A$ and $B$ near infinity, there exist constants
$t_0>0$ and $c>0$  such that
\begin{equation}\label{febbraio106} \overline A(t) \leq A(c t) \quad \hbox{if $t\geq t_0$,} \quad B(t) \leq \overline B(c t) \quad \hbox{if $t\geq
t_0$.}
\end{equation}
From \eqref{febbraio105} and \eqref{febbraio106} one infers that
\begin{align}\label{febbraio107}
\int_\Omega B(|\nabla\bfu|)\,dx & = \int_{\{|\nabla \bfu|<t_0\}}
B(|\nabla\bfu|)\,dx + \int_{\{|\nabla \bfu|\geq t_0\}}
B(|\nabla\bfu|)\,dx \\
\nonumber & \leq  B(t_0) |\Omega| + \int_{\Omega} \overline
B(c|\nabla\bfu|)\,dx
 \leq  B(t_0) |\Omega| + \int _\Omega  \overline A(C c|
\mathcal E^D \bfu |)\, dx
\\
\nonumber & \leq  B(t_0) |\Omega| + \int _{\{C c| \mathcal E^D \bfu
|<t_0\}} \overline A(C c| \mathcal E^D \bfu |)\, dx
\\ \nonumber & \quad
+ \int _{\{C c| \mathcal E^D \bfu |\geq t_0\}} \overline A(C c|
\mathcal E^D \bfu |)\, dx
\\
\nonumber & \leq  (B(t_0) + A(ct_0)) |\Omega| + \int _{\Omega} A(C
c^2 | \mathcal E^D \bfu |)\, dx\,,
\end{align}
namely \eqref{1.3.1} \qed

\medskip
\par\noindent
{\bf Proof of Theorem \ref{thm:main}}. The proof of the equivalence
of (i) and (ii) is completely analogous to that of the corresponding
equivalence in Theorem \ref{thm:main0}, save that Theorem
\ref{sufficiency0} has to be replaced with Theorem
\ref{sufficiency}, and condition (ii) in Corollary
\ref{necessitycor} has to be replaced with condition (iii). The fact
that (iii) implies (ii), and the fact that (i) implies (iii) can be
established along the same lines as in the corresponding
implications in Theorem \ref{thm:main0}. The details are omitted for
brevity. \qed

\medskip
\par\noindent
{\bf Proof of Theorem \ref{thm:plainkorn}}. The derivation of
inequalities \eqref{1.3.plainkorn}
  and \eqref{1.3plainmainkorn} from conditions \eqref{1.1} and
  \eqref{1.2} is the object of \cite[Theorem 3.1]{Ci2} and \cite[Theorem 3.3]{Ci2}, respectively.
Conversely, Theorem \ref{necessity} tells us that inequality
\eqref{1.3.plainkorn}
   implies \eqref{1.1} and
  \eqref{1.2}. Moreover, inequality \eqref{1.3plainmainkorn} implies inequalities  \eqref{1.1} and
  \eqref{1.2} by
  Corollary \ref{necessitycor}, Part (i).
  \qed

\medskip
\par\noindent
{\bf Proof of Theorem \ref{negativenecsuf}}. The validity of
inequality \eqref{march1} under assumptions \eqref{1.1} and
\eqref{1.2} is established in \cite[Theorem 3.1]{BrCi}. We have thus
only to show that \eqref{march1} implies \eqref{1.1} and
\eqref{1.2}. To this purpose, let us introduce negative norms for
single partial derivatives as follows. Given $u\in L^1(\Omega)$, we
set
\begin{align*}
\left\|\tfrac{\partial u}{\partial x_k}\right\|_{W^{-1,A}(\Omega)}
=\sup_{\varphi\in C^\infty_0(\Omega)} \frac{\int_\Omega
u\,\frac{\partial \varphi }{\partial x_k} \,dx}{\|\nabla
\varphi\|_{L^{\widetilde{A}}(\Omega, \rn)}} \quad \hbox{for $k =
1,\dots ,n$}.
\end{align*}
Obviously,
\begin{align}\label{may1}
\left\|\tfrac{\partial u}{\partial x_k}\right\|_{W^{-1,A}(\Omega)}
\leq \|\nabla u\|_{W^{-1,A}(\Omega, \rn)} \quad \hbox{for $k =
1,\dots ,n$}.
\end{align}
On the other hand,
\begin{align}\label{may2}
\|\nabla u\|_{W^{-1,A}(\Omega, \rn)}&=\sup_{\bfvarphi\in
C^\infty_0(\Omega,\rn)} \frac{\int_\Omega u\,\mathrm{div\,}\bfvarphi
\,dx}{\|\nabla \bfvarphi\|_{L^{\widetilde{A}}(\Omega, \rmat)}}
=\,\sup_{\bfvarphi\in C^\infty_0(\Omega,\rn)}
\sum_{k=1}^n\frac{\int_\Omega u\,\frac{\partial \varphi_k}{\partial
x_k} \,dx}{\|\nabla \bfvarphi\|_{L^{\widetilde{A}}(\Omega,
\rmat)}}\\ \nonumber &\leq\,\sup_{\bfvarphi\in
C^\infty_0(\Omega,\rn)} \sum_{k=1}^n\frac{\int_\Omega
u\,\frac{\partial \varphi_k}{\partial x_k} \,dx}{\|\nabla
\varphi_k\|_{L^{\widetilde{A}}(\Omega, \rn)}} \leq\,
\sum_{k=1}^n\sup_{\varphi\in C^\infty_0(\Omega)}\frac{\int_\Omega
u\,\frac{\partial \varphi}{\partial x_k} \,dx}{\|\nabla
\varphi\|_{L^{\widetilde{A}}(\Omega, \rn)}}\\
\nonumber &=\,\sum_{k=1}^n\left\|\tfrac{\partial u}{\partial
x_k}\right\|_{W^{-1,A}(\Omega)},
\end{align}
where $\varphi _k$ denotes the $k$-th component of $\bfvarphi$.
Next, notice   the identity
\begin{equation}\label{may3}
\frac{\partial ^2 v_i }{\partial x_k
\partial x_j}  = \frac{\partial (\mathcal E  \bfv)_{ij} }{\partial
x_k} + \frac{\partial (\mathcal E
 \bfv)_{ik} }{\partial x_j} - \frac{\partial (\mathcal E  \bfv)_{jk}
}{\partial x_i} \end{equation}
 for every weakly differentiable
function $\bfv : \Omega \to \rn$.
\\ Thus, the following chain holds for every
 $\bfu \in W^{1,1}(\Omega, \rn) \cap E^A (\Omega,
\rn)$:
\begin{align}\label{neg1}
\|\nabla \bfu - & (\nabla \bfu )_\Omega \|_{L^B(\Omega , \rmat)}\leq
\,C\,\sum_{i,j=1}^n\left\|\tfrac{\partial u_i}{\partial x_j} -
\big(\tfrac{\partial u_i}{\partial x_j} \big)_\Omega
\right\|_{L^B(\Omega )}\\&\leq \,C\,\sum_{i,j=1}^n\left\|\nabla
\tfrac{\partial u_i}{\partial x_j} \right\|_{W^{-1,A}(\Omega , \rn)}
\leq \,C\,\sum_{i,j,k=1}^n\left\|\tfrac{\partial ^2 u_i }{\partial
x_k
\partial x_j} \right\|_{W^{-1,A}(\Omega)}\nonumber
\\&\leq\,C\,\sum_{i,j,k=1}^n\left(\left\|\tfrac{\partial (\mathcal E\bfu)_{ij}}{\partial x_k}\right\|_{W^{-1,A}(\Omega)}+
\left\|\tfrac{\partial (\mathcal E\bfu)_{ik}}{\partial x_j}\right\|_{W^{-1,A}(\Omega)}+\left\|\tfrac{\partial (\mathcal E\bfu)_{jk}}{\partial x_i}\right\|_{W^{-1,A}(\Omega)}\right)\nonumber\\
&\leq\,C\,\sum_{i,j=1}^n\|\nabla(\mathcal E\bfu)_{ij}\|_{W^{-1,A}(\Omega,\rn)}\leq\,C\,\sum_{i,j=1}^n\|(\mathcal E\bfu)_{ij}-((\mathcal E\bfu)_{ij})_\Omega\|_{L^{A}(\Omega)}\nonumber\\
&\leq \,C\,\|\mathcal E \bfu - (\mathcal E \bfu )_\Omega \|_{L^A(\Omega ,
\rmat)}\nonumber
\end{align}
where the constant $C$ may be different at each occurrence. Note
that the second inequality holds by \eqref{march1}, the third by
\eqref{may2}, the fourth by \eqref{may3}, the fifth by \eqref{may1},
and the sixth by \eqref{negtrivial}.
%
If, in particular, $\bfu \in W^{1,1}_0 (\Omega, \rn)$, then $(\nabla
\bfu )_\Omega = (\mathcal E \bfu )_\Omega =0$, and inequality
\eqref{neg1} implies that
\begin{align}\label{neg2}
\|\nabla \bfu \|_{L^B(\Omega , \rmat)} \leq C  \|\mathcal E
\bfu\|_{L^A(\Omega , \rmat)}
\end{align}
for some constant $C$. The conclusion follows via Theorem
\ref{necessity}, owing to the arbitrariness of $\bfu$. \qed

\medskip
\par\noindent
{\bf Proof of Theorem \ref{corbogov}}. The fact that conditions
\eqref{1.1} and \eqref{1.2} imply inequality \eqref{bogov1} is
proved in \cite[Theorem 3.6]{BrCi}. As far as the converse
implication is concerned, a close inspection of \cite[Inequality
(3.88)]{BrCi} reveals that inequality \eqref{bogov1} implies
inequality \eqref{march1}. The conclusion thus follows from Theorem
\ref{negativenecsuf}. \qed

\bigskip
 \noindent\textbf{Acknowledgments.} 
This research was partly supported by
the  research project of MIUR
(Italian Ministry of Education, University and Research) Prin 2012,
n. 2012TC7588,  ``Elliptic and parabolic partial differential
equations: geometric aspects, related inequalities, and
applications",  by GNAMPA of the Italian INdAM (National
Institute of High Mathematics).

\smallskip
\noindent
The authors wish to thank the referees for their careful reading of  the original version of the paper,  for several valuable suggestions on the presentation, and for pointing out various missing relevant references.


\begin{thebibliography}{[M]}
\bibitem{AcMi} E.Acerbi \& G.Mingione, Regularity results for stationary electro-rheological fluids, \emph{Arch. Rat. Mech. Anal.} {\bf 164} (2002),
213--259.
	\bibitem{AsMa} G. Astarita \& G. Marucci, Principles of non-Newtonian Fluid Mechanics. McGraw-Hill, London,  1974.
\bibitem{Ba} R.Bartnik \& J.Isenberg, The constraint equations. In: P. T. Chrusciel,
H. Friedrich (eds.). The Einstein equations and the large scale behavior of
gravitational fields: 50 years of the Cauchy problem in general relativity,
pp. 1–38, Birkhäuser-Verlag, Basel, Boston, Berlin, 2004.

\bibitem{FB}
M.Bildhauer \& M.Fuchs,  Compact embeddings of the space of
functions with bounded logarithmic deformation, \emph{J.~Math.~Sci.
} \textbf{172} (2011), 165--183.
\bibitem{Bog}
  M.~E.Bogovski{\u\i}, Solutions of some problems of vector analysis, associated with the
  operators ${\rm div}$\ and ${\rm grad}$,
   In {\em Theory of cubature formulas and the application of functional
    analysis to problems of mathematical physics (Russian)}, pages 5--40, \textbf{ 149},
  Akad. Nauk SSSR Sibirsk. Otdel. Inst. Mat., Novosibirsk, 1980.
\bibitem{BrThesis}
D.Breit, \emph{Existence theory for generalized Newtonian fluids},
Postdoctoral thesis, LMU Munich, Department of Mathematics, 2013.

\bibitem{BrCi} D.Breit \& A.Cianchi, Negative Orlicz-Sobolev norms and strongly nonlinear elliptic systems in fluid mechanics, \emph{J. Diff. Eq.}
 \textbf{259}  (2015), 48--83.
\bibitem{BrD} D.Breit \& L.Diening, Sharp conditions for Korn inequalities in Orlicz spaces, \emph{J. Math. Fluid Mech.} \textbf{14}  (2012), 565--573.
\bibitem{BrDF} D.Breit, L.Diening \& M. Fuchs, Solenoidal Lipschitz truncation and applications in fluid mechanics, \emph{J. Diff. Eq.} \textbf{253}
(2012), 1910--1942.
\bibitem{BrSc} D.Breit \&  O.D.Schirra, Korn-type inequalities in Orlicz-Sobolev spaces
involving the trace-free part of the symmetric gradient and
applications to regularity theory, \emph{J. Anal. Appl. (ZAA)}
\textbf{31}  (2012), 335--356.

\bibitem{BreSco94}
S.C.Brenner, \& L.R.Scott, The mathematical theory of finite element
methods, \emph{Texts in Applied Mathematics}, \textbf{15},
Springer-Verlag, New York, (1994) xii+294.


\bibitem{BMM} M.Bul\'i\v cek, M.Majdoub, J.M\'alek, Unsteady flows of fluids with pressure dependent viscosity in unbounded domains,
 \emph{Nonlinear Anal. Real World Appl.} {\bf 11} (2010), 3968--3983.


\bibitem{Cianchisharp} A.Cianchi, A sharp embedding theorem for
Orlicz-Sobolev spaces, \emph{Indiana Univ. Math. J.} \textbf{45
}(1996), 39--65.

\bibitem{Cianchistrong} A.Cianchi, Strong and weak type
inequalities for some classical operators in Orlicz spaces, \emph{J.
London Math. Soc.}  {\bf 60} (1999), 187--202.


 \bibitem{Ci2} A.Cianchi, Korn type inequalities in Orlicz spaces, \emph{J. Funct.
 Anal.}  \textbf{267}  (2014), 2313--2352.
\bibitem{CFM} S.Conti, D.Faraco \& F.Maggi
   A New Approach to Counterexamples to L1 Estimates: Korn's
  Inequality, Geometric Rigidity, and Regularity for Gradients of
  Separately Convex Functions, \emph{Arch. Rat. Mech. Anal.} \textbf{175} (2005), 287--300.
\bibitem{Da} S.Dain, Generalized Korn's inequality and conformal Killing vectors, Calc. Var.
Partial Differential Equations \textbf{25} (2006),   535--540.
\bibitem{DurMus04}
R.G.Dur{\'a}n \& M.A. Muschietti, The {K}orn inequality
  for {J}ones domains, \emph{Electron. J. Differential Equations} \textbf{10} (2004),  10 pp. (electronic).
\bibitem{DRS} L.Diening, M.R\r{u}\v{z}i\v{c}ka \&
  K.Schumacher, A Decomposition technique for John
  domains, \emph{Ann. Acad. Scientiarum Fennicae} \textbf{35}  (2009), 87--114.
\bibitem{fei3} E. Feireisl, \emph{Dynamics of Compressible Flow}, 
Oxford University Press, Oxford, 2004.
\bibitem{FeMa} E. Feireisl, X. Liao \& J. M\'alek, Global weak solutions to a class of non-Newtonian compressible fluids, \emph{Math. Meth. Appl. Sci.} \textbf{38} (2015), 3482--3494.
\bibitem{F} E. Feireisl \& A. Novotn\'{y}, \emph{Singular limits in thermodynamics of viscous fluids.} Birkh\"auser-Verlag, Basel, 2009.
\bibitem{E} H. J. Eyring (1936): Viscosity, plasticity, and diffusion as example of absolute reaction rates. J. Chemical Physics  4, 283-291.
\bibitem{FZ} D.Faraco \& X.Zhong, Geometric rigidity of conformal matrices. \emph{Ann. Sc. Norm. Super. Pisa Cl. Sci.}  \textbf{4} (2005), 557--585.
    \bibitem{FrS} J.Frehse \& G.Seregin,  Regularity of solutions
    to variational problems of the deformation theory of plasticity with logarithmic hardening,
    \emph{ Proc. St. Petersburg Math. Soc.} \textbf{5}, 184--222;
     English Translation: \emph{Amer. Math. Soc. Transl.} II \textbf{193} (1998/1999), 127--152
\bibitem{FJM1} G.Friesecke, R.D.James \& S.M\"uller, A theorem on geometric rigidity and
the derivation of nonlinear plate theory from three-dimensional elasticity, \emph{Comm. Pure
Appl. Math.}  \textbf{55} (2002),   1461--1506.
\bibitem{FJM2} G.Friesecke, R.D.James \& S.M\"uller, A hierarchy of plate models derived
	from nonlinear elasticity by gamma-convergence, \emph{Arch. Ration. Mech. Anal.} \textbf{180} (2006),
183--236.
\bibitem{Fu1} M.Fuchs, On stationary incompressible
  Norton fluids and some extensions of Korn's
  inequality, \emph{Zeitschr. Anal. Anwendungen} \textbf{13} (1994), 191--197.
\bibitem{Fu2} M.Fuchs,  Korn
  inequalities in Orlicz spaces, \emph{Irish Math. Soc. Bull. } \textbf{65} (2010), 5--9.
\bibitem{FR}
 M.Fuchs \& S.Repin, Some Poincar\'{e}-type inequalities for
 functions of bounded deformation involving the deviatoric part of
 the symmetric gradient, \emph{Zap. Nauchn. sem. St.-Petersburg
 Odtel. Math. Inst. Steklov (POMI)} \textbf{385} (2010), 224--234.
\bibitem{FuSc}
M.Fuchs \& O.Schirra, An application of a new coercive inequality to
variational problems studied in general relativity and in Cosserat
elasticity giving the smoothness of minimizers, {\it Arch.~Math.}
\textbf{93} (2009), 587--596.
\bibitem{FuS} M.Fuchs \& G.Seregin,
 Variational methods for problems from plasticity theory and for generalized Newtonian fluids. Lecture Notes in Mathematics Vol. 1749,
 Springer Verlag, Berlin-Heidelberg-New York,  2000.
\bibitem{Go1} J.Gobert, Une in{\'e}quation
  fondamentale de la th{\'e}orie de l'{\'e}lasticit{\'e},
  \emph{Bull. Soc. Roy. Sci. Liege} \textbf{3-4}  (1962), 182--191.
\bibitem{Go2} J.Gobert, Sur une in{\'e}galit{\'e} de
  coercivit{\'e}, \emph{J. Math. Anal. Appl.} \textbf{36}  (1971), 518--528.
\bibitem{HRJ} T.A.Hassan, V.K.Rangari \& S.Jeelani, Synthesis, processing and characterization of shear thickening fluid (STF) impregnated fabric
composites, \emph{Materials Science and Engineering: A} \textbf{527}
(2010), 2892--2899.
\bibitem{Je} J.Jeong, H.Ram\'ezani, I.M\"unch \& P.Neff, A numerical study for linear isotropic Cosserat elasticity with
conformally invariant curvature, \emph{Z. Angew. Math. Mech.}
\textbf{89} (2009), 552--569.
\bibitem{KMS} B.Kirchheim, S.M{\"u}ller \& V.\v{S}v{\'e}rak, Studying
  nonlinear pde by geometry in matrix space. Geometric analysis and
  nonlinear partial differential equations (S.Hildebrandt, H.Karcher
  eds.), Springer  (2003), 347-395.
\bibitem{KK} V.Kokilashvili \& M.Krbec,
 \lq\lq Weighted inequalities in Lorentz and Orlicz spaces",  World Scientific Publishing, River Edge, NJ, 1991.
\bibitem{KO} V.A.Kondratiev \& O.A.Oleinik, On Korn's inequalities,
C. R. Acad. Sci. Paris Ser. I \textbf{308} (1989), 483--487.
\bibitem{Korn} A.Korn, \"Uber einige Ungleichungen, welche in der Theorie der
elastischen und elektrischen Schwingungen eine Rolle spielen, in:
Classe des Sciences Math\'ematiques et Naturels (9, Novembre),
 \emph{Bull. Internat. Acad. des Sci. Cracovie}  (1909), 705--724.
\bibitem{LMu} M.Lewicka \& S.M\"uller, The uniform Korn-Poincar\'e inequality in thin domains,
\emph{Ann. Inst. H. Poincar\'e Anal. Non Lin\'eaire} \textbf{28} (2011), 443--469.
   \bibitem{LM} K.de Leeuw \& H.Mirkil, A priori estimates for differential operators in $L_\infty$ norm, \emph{Illinois J. Math.} \textbf{8}  (1964), 112--124.
\bibitem{Li2} P.\,L. Lions, \emph{Mathematical topics in fluid mechanics. Vol. 2. Compressible models.} Oxford Science Publications, The Clarendon Press, Oxford University Press, New York,  1998.
\bibitem{MaRa} J. M\'alek \& K.R. Rajagopal,  Compressible generalized Newtonian fluids,  \emph{Zeit. Angew. Math.   Physik (ZAMP)} \textbf{61} (2010), 1097--1110.
\bibitem{M} S. M{\"u}ller, Variational models for microstructure
  and phase transitions, In: \emph{Calculus of variations and geometric
  evolution problems} (F. Bethuel, ds.), Springer Lecture Notes in
  Math. 1713, Springer, Berlin  (1999), 85-210.
\bibitem{Mik} S.G.Mikhlin,  \emph{Multidimensional singular integrals and integral equations}, Pergamon press, Oxford,  1965.
\bibitem{MM} P.~P.Mosolov \& V.~P.Mjasnikov, On the correctness
  of boundary value problems in the mechanics of continuous
  media, \emph{Math. USSR Sbornik} \textbf{17}  (1972), 257--267.
\bibitem{Ne} J.Ne\v{c}as, Sur les normes {\'e}quivalentes dans
  $W_k^p(\Omega)$ et sur la coecivit{\'e} des formes formellement
  positives, in S{\'e}minaire Equations aus D{\'e}riv{\'e}es Partielles, Les
  Presses de l'Universit{\'e} de Montr{\'e}al  (1966), 102--128.
\bibitem{Ne} P.Neff, D.Pauly \& K.-J.Witsch, Poincare meets Korn via Maxwell: extending Korn's first inequality to incompatible tensor fields, \emph{J. Diff. Equat.} \textbf{ 258}
(2015), 1267--1302.
\bibitem{NeJe} P.Neff \&
J.Jeong, A new paradigm: the linear isotropic Cosserat model
with conformally invariant curvature energy, \emph{Z. Angew. Math. Mech.}
\textbf{89} (2009), 107--122.
\bibitem{NeJe2} P.Neff, J.Jeong \& A. Fischle, Stable identification of linear isotropic Cosserat parameters: bounded stiffness in bending and torsion implies conformal invariance of curvature, \emph{ Acta Mechanica} \textbf{ 211} (2010), 237--249.
\bibitem{Or} D.Ornstein, A non-inequality for
  differential operators in the $L_1$ norm, \emph{Arch. Rat.
  Mech. Anal.} \textbf{11}  (1964), 40-49.
\bibitem{RR1}  M.M.Rao \& Z.D.Ren, \emph{Theory of Orlicz spaces}, Marcel Dekker Inc.,
 New York, 1991.
 \bibitem{RR2}  M.M.Rao \& Z.D.Ren, \emph{Applications of Orlicz spaces}, Marcel Dekker Inc.,
 New York, 2002.

\bibitem{Resh} Yu.G.Reshetnyak, Estimates for certain
differential operators with finite dimensional kernel,
\emph{Sibirskii Math. Zh.} \textbf{2}  (1970), 414--418.
\bibitem{Resh2} Yu.G.Reshetnyak, \emph{Stability theorems in geometry and analysis}, Kluwer Academic
Publishers Group, Dordrecht, 1994.

\bibitem{Sc2}
O.Schirra,  New Korn-type inequalities and regularity of solutions
to linear elliptic systems and anisotropic variational problems
involving the trace-free part of the symmetric gradient, {\it
Calc.~Var.} \textbf{43} (2012), 147 -- 172.

\bibitem{SMB} A.Srivastavaa, A.Majumdara \& B.S.Butolaa,  Improving the Impact Resistance of Textile Structures by
using Shear Thickening Fluids: A Review, \emph{Critical Reviews in
Solid State and Materials Sciences}  \textbf{37}  (2012), 115--129.

\bibitem{Te} R. Temam (1985): Mathematical problems in
  plasticity. Gauthier Villars, Paris.



\bibitem{W} A.Wr\'oblewska, Steady flow of non-Newtonian fluids--monotonicity methods in generalized Orlicz
spaces, \emph{Nonlinear Anal.}   \textbf{72}  (2010), 4136--4147.


\end{thebibliography}
\end{document}